\documentclass[12pt]{article}
\usepackage{amsmath,amsthm,amssymb}
\usepackage{mathrsfs}
\usepackage{color}
\DeclareMathOperator{\qdim}{qdim}

\DeclareMathOperator{\Gal}{Gal}
\DeclareMathOperator{\ch}{ch}
\DeclareMathOperator{\tr}{tr}

\DeclareMathOperator{\glob}{glob}
\begin{document}
\input amssym.def
\setcounter{equation}{0}
\newcommand{\wt}{{\rm wt}}
\newcommand{\spa}{\mbox{span}}
\newcommand{\Res}{\mbox{Res}}
\newcommand{\End}{\mbox{End}}
\newcommand{\Ind}{\mbox{Ind}}
\newcommand{\Hom}{\mbox{Hom}}
\newcommand{\Mod}{\mbox{Mod}}
\newcommand{\m}{\mbox{mod}\ }
\renewcommand{\theequation}{\thesection.\arabic{equation}}
\numberwithin{equation}{section}

\def \End{{\rm End}}
\def \Aut{{\rm Aut}}
\def \Z{\mathbb Z}
\def \H{\mathbb H}
\def \M{\Bbb M}
\def \C{\mathbb C}
\def \R{\mathbb R}
\def \Q{\mathbb Q}
\def \N{\mathbb N}
\def \ann{{\rm Ann}}
\def \<{\langle}
\def \o{\omega}
\def \O{\Omega}
\def \Or{\cal O}
\def \M{{\cal M}}
\def \1t{\frac{1}{T}}
\def \>{\rangle}
\def \t{\tau }
\def \a{\alpha }
\def \e{\epsilon }
\def \l{\lambda }
\def \L{\Lambda }
\def \G{\bar{G}}
\def \g{\gamma}
\def \b{\beta }
\def \om{\omega }
\def \o{\omega }
\def \ot{\otimes}
\def \cg{\chi_g}
\def \ag{\alpha_g}
\def \ah{\alpha_h}
\def \ph{\psi_h}
\def \S{\cal S}
\def \nor{\vartriangleleft}
\def \V{V^{\natural}}
\def \1{{\bf 1}}
\def \be{\begin{equation}\label}
\def \ee{\end{equation}}
\def \qed{\mbox{ $\square$}}
\def \pf {\noindent {\bf Proof:} \,}
\def \bl{\begin{lem}\label}
\def \el{\end{lem}}
\def \ba{\begin{array}}
\def \ea{\end{array}}
\def \bt{\begin{thm}\label}
\def \et{\end{thm}}
\def \br{\begin{rem}\label}
\def \er{\end{rem}}
\def \ed{\end{de}}
\def \bp{\begin{prop}\label}
\def \ep{\end{prop}}
\def \p{\phi}
\def \d{\delta}
\def \irr{{\rm Irr}}
\def\wt{{\rm wt}}
\def\res{{\rm Res}}
\def\CC{{\cal C}}
\def\dim{{\rm dim}}
\def\hom{\rm Hom}
\def\id{\rm id}

\def \voa{vertex operator algebra\ }
\def \vosa{vertex operator superalgebra\ }
\def \vosas{vertex operator superalgebras\ }
\def \voas{vertex operator algebras}
\def \v{vertex operator algebra\ }
\def\AGNV{A_{g,n}(V)}							
\def\OGNV{O_{g,n}(V)}
\def\VRAST{V^{r\ast}}	

\newtheorem{th1}{Theorem}
\newtheorem{ree}[th1]{Remark}
\newtheorem{thm}{Theorem}[section]
\newtheorem{prop}[thm]{Proposition}
\newtheorem{coro}[thm]{Corollary}
\newtheorem{lem}[thm]{Lemma}
\newtheorem{rem}[thm]{Remark}
\newtheorem{de}[thm]{Definition}
\newtheorem{hy}[thm]{Hypothesis}
\newtheorem{conj}[thm]{Conjecture}
\newtheorem{ex}[thm]{Example}

\begin{center}
{\Large {\bf Super orbifold theory}}\\
\vspace{0.5cm}

Chongying Dong\footnote
{Supported by NSFC grant 11871351}\\
Department of Mathematics, University of
California, Santa Cruz, CA 95064 USA \\
Li Ren\footnote{Supported by NSFC grants 11301356 and 11671277}\\
 School of Mathematics,  Sichuan University,
Chengdu 610064 China\\
Meiling Yang\\
School of Mathematics,  Sichuan University,
Chengdu 610064 China

\end{center}

\begin{abstract}
Let $V$ be a simple vertex operator superalgebra and $G$ a finite automorphism group of $V$ contains a canonical automorphism $\sigma$ such that $V^G$ is regular. 
It is proved that every irreducible $V^G$-module occurs in an irreducible $g$-twisted $V$-module for some $g\in G.$
 Moreover, the quantum dimensions of irreducible $V^G$-modules are determined, a global dimension formula for $V$ in terms of twisted modules is obtained and a super quantum Galois theory is established. In addition the $S$-matrix of 
 $V^G$ is computed.
\end{abstract}

\section{Introduction}
This paper is a continuation of our study on orbifold theory for a vertex operator superagebra for an arbitrary finite automorphism group. A vertex operator superalgebra $V=V_{\overline{0}} \oplus V_{\overline{1}}$ in this paper is $\frac{1}{2}\Z$-graded with $V_{\overline{i}}=\frac{i}{2}+\Z$ for $i=0,1.$ Then $V$ has a canonical automorphism $\sigma$ such that
$\sigma =(-1)^i.$ Let  $G$ be a  finite automorphism group of $V$ containing the canonical automorphism $\sigma$.
Then  the $G$-fixed points $V^G$ in $V$  is a vertex operator algebra.  The main purpose of this paper to understand 
the representations of $V^G$ in terms of $G$ and twisted $V$ modules.  If $G$ is generated by $\sigma$, 
the $V^G$-module category and its connection with 16 fold way conjecture have been investigated extensively in \cite{DNR1}.

A systematic study of irreducible $V^G$-modules in an arbitrary vertex operator algebra $V$ started in \cite{DM1,  DLM2}. It turms out that the actions of $G$ and $V^G$ on $V$ form a dual pair in the sese of \cite{Ho}.  This duality result was generalized  in \cite{DY, MT} to include the twisted modules. A well known conjuecture in orbifold theory says that if a vertex operator algebra $V$ is rational, then $V^G$ is rational and any irreducible $V^G$-module occurs in an irreducible $g$-twisted $V$-module for some $g\in G.$ 
It is  proved in \cite{CM} that  $V^G$ is regular (rational and $C_2$-cofinite)  if  $V$ is regular and $G$ is cyclic.  Recently,  it is established in \cite{DRX1} that if $V^G$ is regualr then all the irreducible $V^G$-module  occur in irreducible twisted $V$-modules. In particular the conjecture holds  for solvable group.

In this paper  we extend the results  for a vertex operator algebra \cite{DRX1} to a \vosa for  any finite automorphism group $G$ which contains $\sigma$ under the assumption that  $V^G$ is regular and the weight of any irreducible $g$-twisted $V$ module except $V$ itself is positive. In particular, we classify the irreducible
$V^G$-modues, determin the quantum dimensions of irreducible $V^G$-modules, obtain a global dimension formula for $V$ in terms of twisted modules and establish a super quantum Galois theory. 
We also compute the $S$-matrix of $V^G.$ The modulaity of trace functions in super orbifold theory in \cite{DZ1} plays an essential role in proving the main results.

Classification of irreducible $V^G$-modules consists of two parts: 1) Every irreducible $V^G$-module appers in an irreducible $g$-twisted $V$-module for some $g\in G,$ 2) Determine when the two irreducible $V^G$-modules 
appearing in twisted $V$-modules  are isomorphic.  1) is proved by using the modular invariance of trace functions 
in orbifold theory \cite{DZ1} (also see \cite{Z,DLM7}) and nonvanishing property of  the entries $S_{V^G, W}$ 
for any irreducible $V^G$-module $W$ \cite{Hu}. One can write 
the trace function $\tau^{-\wt[v]}Z_{V^G}(v,\frac{-1}{\tau})$ as a linearly combination of linearly 
indepedent trace functions $Z_W(v,\tau)$ on $V^G\times \H$ where $W$ are the irreducible $V^G$-modules occuring in twisted $V$-modules and $\H$ is the upper half plane, and conclude that these are all the irreducible $V^G$-modules. Identifing irreducible $V^G$-modules is more involved. Let $ {\cal S}=\cup_{g\in G} {\cal S}(g)$ where $ {\cal S}(g)$ is the set of inequivalent of irreducible $g$-twisted $V$-modules. Then $G$ acts on $ {\cal S}$  such that $M\circ k=M$
as a vector space and $Y_{M\circ k}(v,z)=Y_M(kv,z)$ where $(M,Y_M)$ is an irreducible $g$-twisted $V$-module, $k\in G.$ Then $M\circ k$ is an irrducible $k^{-1}gk$-twisted $V$-module and the stablizer $G_M=\{k\in G|M\circ k\cong M\}$ 
acts on $M$ projectively.  So one can decompose $M$ into $\oplus_{\l\in \Lambda_M}W_\l\otimes M_\l$ where
$W_{\l}$ are irreducible projective $G_M$-modules and $M_{\l}$ are irreducibe $V^{G_M}$-modules.
Let  ${\cal S}=\cup_{j\in J}M^j\circ G$ is a union of disjoint orbits. Then $\{M^j_{\l}|j\in J,\l\in \Lambda_{M^j}\}$ classifies the irreducible $V^G$-modules. This result is similar to the classification of irreducible modules in  orbifold theory for vertex operator algebras  \cite{DRX1}, but the proof is more complex and subtle.  While classification of irreducible $V^G$-modules  uses irreducible twisted $V$-modules, the modular invariance result in \cite{DZ1} deals with only irreducible twisted super modules. Let $\overline{G}$ be the restriction of $G$ to vertex operator algebra $V_{\bar 0}.$ Then $V^G=V^{\overline{G}}_{\bar 0}$ and $\overline{G}$ acts on the set $\overline{\cal S}$ of inequivalent irreducible $\overline{G}$-twisted $V_{\bar 0}$-modules for $g\in G$ in the same way. Another effort is to find relations between the projective respresentations of a stablizer $G_M$ for $M\in {\cal S}$ and the projective representations of the stablizers $\overline{G}_{\bar M}$ of relavant objects $\bar M$
in  $\overline{\cal S}.$ 

The quantum dimensions $\qdim_V M$ of an irreducible $g$-twisted $V$-module $V$ is defined to be the limit of $\frac{\ch_q M}{\ch_q V}$ as $q$ goes to 1 as in \cite{DJX} where the quantum dimension of a module for vertex operator algebras was first defined and studied.  Roughly speaking,  the quantum dimension $\qdim_V M$ is equal to $\frac{\dim M}{\dim V}$ of type $\frac{\infty}{\infty}.$ Using the $S$-matrix from \cite{DZ1} we give an explict formula 
for the $\qdim_V M.$ In particular,  $\qdim_V M$ exists.  The relation between the quantum dimension of an irreducible $g$-twisted $V$-module $M$ and the quantum dimension of an irreducible $V^G$-submodule $M_{\l}$ 
of $M$ is given by
$$\qdim_{V^G}M_{\l}=[G:G_M]\dim W_{\l}\qdim_VM.$$ 
In particular $\qdim_{V^G}V_{\l}=\dim W_{\l}$ for any irredicible character $\l$ of $G.$  
This leads to a super quantum Galois correspondence in the sense that $H\mapsto V^G$ yields a bijection between subgroups $H$ of $G$ and 
vertex operator super subalgebras of $V$ contaning $V^G$ satisfying   $[V:V^H]=o(H)$ and $[V^H:V^G]=[G:H]$
where $[V:U]=\qdim_UV$ if $U$ is a conformal subalgebra of $V.$ Moreover, the Galois group $Gal(V/V^G)=G.$
The interesting part is to show that $H\mapsto V^H$ is onto. If $U$ is a vertex operator subalgebra of $V_{\bar 0}$ it is easy
to show that there exists $H<G$ such that $U=V^H$ by applying the Galois correpondence  in \cite{DM1, HMT,DJX} to vertex operator algebra $V_{\bar 0}.$ If $U$ is not a subagebra of $V_{\bar 0}$ the proof is nontrivial. 

Analogous to the dimension of a finite dimensional semisimple associaitve algebra, the global dimension $\glob(V)$ of a \vosa $V$ is defined to be the sum of $(\qdim_VM)^2$ over inequivalent irreducilbe $V$-modules $M$ \cite{DJX}.  As in \cite{DRX1} we show that $\glob(V)$ is also the sum  of $(\qdim_VM)^2$ over inequivalent irreducilbe $g$-twisted $V$-modules $M$ for any automorphism $g$ of finite order, and that $\glob(V^G)=o(G)^2\glob(V)$ for any finite automorphism group $G.$ It is known that the $V^G$-module category $\CC_{V^G}$ is a modular tesnor category  \cite{Hu} and $V$ is a commutative algebra in the $V^G$-module category \cite{HKL, CKM}. There is a fusion category $(\CC_{V^G})_V$  of left  $V$-modules in 
category  $\CC_{V^G}$ \cite{KO}. If $V$ is a vertex operator algebra,   the simple objects in $(\CC_{V^G})_V$ are exactly the irreducible $g$-twisted $V$-modules for $g\in G$ \cite{DLXY}. In this case the quantum dimensions are, in fact,  the Frobenius-Perron dimensions in the fusion category $(\CC_{V^G})_V.$  We expect that this is also true if $V$ is a vertex operator superalgebra. 

It is worthy to discuss more on super twisted modules.   A $g$-twisted $V$-module is called a super twisted $V$-module if $M\circ\sigma\cong M.$ An irreducible super twisted module is either an irreducibe twisted module or a direct sum of two inequivalent irreducible twisted modules.  We now explian even when $V$ is holomorphic in the sense that  $V$ is the only irreducible module for itself up to isomorphism,  the orbifold theories  for vertex operator algebra and 
vertex operator superalgebra are very different. It is known from \cite{DLM7} that
if $V$ is a holomorphic vertex operator algebra then for any finite order automorphism $g$ there is a unique irreducible $g$-twisted $V$-module $V(g)$ which has quantum dimension 1 and is a simple current. Moreover, the $V^G$-module
category $\CC_{V^G}$ is braided equivalent to the module category of twisted Drinfeld double $D^{\a}[G]$ for some
$\alpha\in H^3(G,S^1)$ \cite{DPR, K, DNR2}. Although a holomorphic \vosa $V$ has a unique irreducible super $g$-twisted module \cite{DZ1},  it can have two inequivalent irreducible $g$-twisted modules which have quantum dimensions $\frac{1}{\sqrt{2}}$ as the global dimension of $V$ is 1. The $V^G$-module category in this case has not been understood except for $G=\Z_2$ \cite{Ki, DNR1}.

The paper is organized as follows. In Section 2, we review various notions of twisted modules, $g$-rationality following \cite{DLM4, DZ2}. We recall  the twisted associative algebra $A_{g,n}(V)$ in Section 3 for $0\leq n\in\frac{1}{T'}\Z$ and relevant results from \cite{P} for the purpose of classification of irreducible $V^G$-modues later. Section 4 is a review of the modular invariance of trace functions from \cite{DZ1}. The transformation of these functions by the $S=\left(\begin{array}{cc} 0 &-1\\ 1 &0\end{array}\right)$ plays an essential role in the computation of the quantum dimensions and the global dimensions.
In Section 5, we first classify the irreducible $V^G$-modules appearing in twisted $V$-modules and then prove 
that these are all the irreducible $V^G$-modules.
Section 6 is devoted to the study of quantum dimensions. The quantum dimension of any irreducible $g$-twisted module $M$ is computed by using the $S$-matrix. The quantum dimension of any irreducible $V^G$-submodule of $M$  is obtained then by using the quantum dimension of $M.$ Section 7 is about the global dimension.  The global dimensions of $V$ and $V^G$ are given explicitly in terms of the $S$-matrix from \cite{DZ1}. In Section 8, we compute the $S$-matrix of $V^G$ following \cite{DRX2}.

\section{Basics}
We review  various notions of twisted modules for a \vosa following \cite{DZ1, DZ2}   in this section. 
We also discuss important concepts such as  rationality, regularity, and $C_2$-cofiniteness from \cite{Z,DLM3} . 

A super vector space  is a $\Bbb Z_{2}$-graded vector space
 $V=V_{\bar{0}}\oplus V_{\bar{1}}$.  The vectors  in $V_{\bar{0}}$
(resp. $V_{\bar{1}}$) are called even (resp. odd). Let $\tilde{v}$
be $0$ if $v\in V_{\bar{0}}$, and $1$ if  $v\in V_{\bar{1}}$.

\begin{de} \label{SVOA}
 A  vertex operator superalgebra is a
$\frac{1}{2}\Bbb Z$-graded super vector space
\begin{equation*}
V=\bigoplus_{n\in{ \frac{1}{2}\Bbb Z}}V_n= V_{\bar{0}}\oplus V_{\bar{1}}
\end{equation*}
with  $V_{\bar{0}}=\oplus_{n\in\Z}V_n$ and
$V_{\bar{1}}=\oplus_{n\in\frac{1}{2}+\Z}V_n$
satisfying $\dim V_{n}< \infty$ for all $n$ and $V_n=0$ if $n$ is sufficiently
small.  The space $V$ is equipped with a linear map
 \begin{align*}
& V \to (\mbox{End}\,V)[[z,z^{-1}]] ,\\
& v\mapsto Y(v,z)=\sum_{n\in{\Z}}v_nz^{-n-1}\ \ \ \  (v_n\in
\mbox{End}\,V)\nonumber
\end{align*}
and with two distinguished vectors ${\bf 1}\in V_0,$ $\omega\in
V_2$ satisfying the following conditions for $u, v \in V,$ and $m,n\in\Z:$
\begin{align*} \label{0a4}
& u_nv=0\ \ \ \ \ {\rm for}\ \  n\ \ {\rm sufficiently\ large};  \\
& Y({\bf 1},z)=Id_{V};  \\
& Y(v,z){\bf 1}\in V[[z]]\ \ \ {\rm and}\ \ \ \lim_{z\to
0}Y(v,z){\bf 1}=v;\\
& [L(m),L(n)]=(m-n)L(m+n)+\frac{1}{12}(m^3-m)\delta_{m+n,0}c ;\\
& \frac{d}{dz}Y(v,z)=Y(L(-1)v,z);\\
& L(0)|_{V_n}=n
\end{align*}
where $L(m)=\o_{ m+1}, $ that is,
$$Y(\o,z)=\sum_{n\in\Z}L(n)z^{-n-2};$$
and  {\em Jacobi identity} holds:
\begin{equation*}\label{2.8}
\begin{array}{c}
\displaystyle{z^{-1}_0\delta\left(\frac{z_1-z_2}{z_0}\right)
Y(u,z_1)Y(v,z_2)-(-1)^{\tilde{u}\tilde {v}}z^{-1}_0\delta\left(\frac{z_2-z_1}{-z_0}\right)
Y(v,z_2)Y(u,z_1)}\\
\displaystyle{=z_2^{-1}\delta
\left(\frac{z_1-z_0}{z_2}\right)
Y(Y(u,z_0)v,z_2)}
\end{array}
\end{equation*}
where $\delta(z)=\sum_{n\in\Z}z^n$ and
all binomial expressions (here and below) are to be expanded in nonnegative
integral powers of the second variable.
\end{de}
Such a vertex operator superalgebra may be denoted by  $V=(V,Y,{\bf 1},\omega).$
If $V_{\bar 1}=0,$ $V$ is a vertex operator algebra $V$ \cite{B,FLM}.

A invertible linear transformation $g$ of \vosa  $V$ is called an automorphism if $g\1=\1,$ $g\omega=\omega$ and
$gY(u,z)g^{-1}=Y(gu,z)$ for all $u\in V.$ It is obvious that $g$ preseres  both $V_{\bar 0}$ and $V_{\bar 1}.$
There is a canonical automorphism $\sigma$ of $V$ satisfying  $\sigma|_{V_{\bar i}}=(-1)^i$ for $i=0,1.$
Clearly, $\sigma$ is a central element of $\Aut(V).$ 

Now fix an automorphism $g$ of $V$ of order $T<\infty$. Let $o(\sigma g)=T'.$ 
Then  $V$ decomposes  into eigenspaces for $g$ and $\sigma :$
\begin{equation*}\label{g2.1}
V=\bigoplus_{r\in \Z/T\Z}V^r ,  \
V=\bigoplus_{r\in \Z/T'\Z}V^{r*}
\end{equation*}
where $V^r=\{v\in V|gv=e^{-2\pi ir/T}v\}$, and $V^{r*}=\{v\in V| \sigma g v=e^{-2\pi ir/T'}v\}.$
We use $r$ to denote both
an integer between $0$ and $T-1$ (resp. $T'-1$ ) and its residue class modulo $T$  (resp. $T'$ ) in this
situation.

\begin{de} \label{weak}
A {\em weak $g$-twisted $V$-module} $M$ is a vector space equipped
with a linear map
\begin{equation*}
\begin{split}
Y_M: V&\to (\End\,M)[[z^{1/T},z^{-1/T}]]\\
v&\mapsto\displaystyle{ Y_M(v,z)=\sum_{n\in\frac{1}{T}\Z}v_nz^{-n-1}\ \ \ (v_n\in
\End\,M)},
\end{split}
\end{equation*}
which satisfies the following:  for all $0\leq r\leq T-1,$ $u\in V^r$, $v\in V,$
$w\in M$,
\begin{eqnarray*}
& &Y_M(u,z)=\sum_{n\in \frac{r}{T}+\Z}u_nz^{-n-1} \label{1/2},\\
& &u_lw=0~~~
\mbox{for}~~~ l\gg 0,\label{vlw0}\\
& &Y_M({\mathbf 1},z)=Id_M,\label{vacuum}
\end{eqnarray*}
 \begin{equation*}\label{jacobi}
\begin{array}{c}
\displaystyle{z^{-1}_0\delta\left(\frac{z_1-z_2}{z_0}\right)
Y_M(u,z_1)Y_M(v,z_2)-(-1)^{\tilde{u}\tilde {v}}z^{-1}_0\delta\left(\frac{z_2-z_1}{-z_0}\right)
Y_M(v,z_2)Y_M(u,z_1)}\\
\displaystyle{=z_2^{-1}\left(\frac{z_1-z_0}{z_2}\right)^{-r/T}
\delta\left(\frac{z_1-z_0}{z_2}\right)
Y_M(Y(u,z_0)v,z_2)},
\end{array}
\end{equation*}
\end{de}


\begin{de}\label{ordinary}
A $g$-{\em twisted $V$-module} is
a $\C$-graded weak $g$-twisted $V$-module
\begin{equation*}
M=\bigoplus_{\lambda \in{\C}}M_{\lambda}
\end{equation*}
where $M_{\l}=\{w\in M|L(0)w=\l w\}$ and $L(0)$ is the component operator of $Y(\omega,z)=\sum_{n\in \Z}L(n)z^{-n-2}.$ We also require that
$\dim M_{\l}$ is finite and for fixed $\l,$ $M_{\frac{n}{T'}+\l}=0$
for all small enough integers $n.$ If $w\in M_{\l}$ we refer to $\l$ as the {\em weight} of
$w$ and write $\l=\wt w.$
\end{de}

We use $\Z_+$ to denote the set of nonnegative integers.
\begin{de}\label{admissible}
 An {\em admissible} $g$-twisted $V$-module
is a  $\frac1{T'}{\Z}_{+}$-graded weak $g$-twisted $V$-module
\begin{equation*}
M=\bigoplus_{n\in\frac{1}{T'}\Z_+}M(n)
\end{equation*}
such that $
v_mM(n)\subseteq M(n+\wt v-m-1)$
for homogeneous $v\in V,$ $m\in \frac{1}{T}{\Z},n\in \frac{1}{T'}{\Z}.$
\ed

In the case that $g=Id_V$  we have the notions of  weak, ordinary and admissible $V$-modules, respectively.

Let  $M=\bigoplus_{n\in \frac{1}{T'}\Z_+}M(n)$ be  an admissible $g$-twisted $V$-module. Then
\begin{equation*}
M'=\bigoplus_{n\in \frac{1}{T'}\Z_+}M(n)^{*}
\end{equation*}
 is an admissible $g^{-1}$-twisted $V$-module 
where $M(n)^*=\Hom_{\C}(M(n),\C)$ and the vertex operator
$Y_{M'}(a,z)$ is defined for $a\in V$ via
\begin{eqnarray*}
\langle Y_{M'}(a,z)f,u\rangle= \langle f,Y_M(e^{zL(1)}(e^{\pi i}z^{-2})^{L(0)}a,z^{-1})u\rangle
\end{eqnarray*}
where $\langle f,u\rangle=f(u)$ is the natural paring $M'\times M\to \C$ \cite{FHL, X}.
Similarly, we can define the contragredient module $M'$ for a $g$-twisted $V$-module $M.$ In this case,
$M'$ is a $g^{-1}$-twisted $V$-module. Moreover, $M$ is irreducible if and only if $M'$ is irreducible.
$V$ is called selfdual if $V'$ and $V$ are isomorphic $V$-modules.

\begin{de} (1) A \vosa $V$ is called $g$-rational, if the  admissible $g$-twisted module category is semisimple. $V$ is called rational if $V$ is $1$-rational.

(2) A \vosa $V$ is $C_2$-cofinite if $V/C_2(V)$ is finite dimensional, where $C_2(V)=\langle v_{-2}u|v,u\in V\rangle.$
\end{de}

\begin{thm}\label{grational}{\cite{DZ1, DZ2}} Asuume that $V$ is $g$-rational. Then

(1) Every irreducible admissible $g$-twisted $V$-module $M$ is a $g$-twisted $V$-module. Moreover, there exists a number $\l \in \mathbb{C}$ such that  $M=\oplus_{n\in \frac{1}{T'}\mathbb{Z_+}}M_{\l +n}$ where $M_{\lambda}\neq 0.$ The $\l$ is called the conformal weight of $M;$

(2) There are only finitely many irreducible $g$-twisted $V$-modules up to isomorphism.

(3) If $V$ is also $C_2$-cofinite and $g^i$-rational for all $i\geq 0$ then the central charge $c$ and the conformal weight $\l$ of any irreducible $g$-twisted $V$-module $M$ are rational numbers.
\end{thm}

\begin{de}
A \vosa $V$ is called regular if every weak $V$-module is a direct sum of irreducible $V$-modules.
\end{de}

A \vosa $V=\oplus_{n\in \Z}V_n$  is said to be of CFT type if $V_n=0$ for negative
$n$ and $V_0=\C {\bf 1}.$  If  $V$ is a vertex operator algebra of CFT type, then regularity is equivalent to rationality and $C_2$-cofiniteness \cite{DLM3, L,ABD, DYu, HA}.

\section{Associative algebras  $\AGNV$ and related results }

We  review  associative algebras $\AGNV$ for $n\in \frac{1}{T'}\Z_+$ studied in \cite{P} (also see \cite{DLM6, DZ2})
where $T'$ is the order of $g\sigma$ as before.  We also present some basic results related to  $\AGNV.$
  If $n=0$ we obtain the associative algebra $A_g(V)$ constructed in \cite{DZ2}. With $g=1$ we recover associative algebra $A_n(V)$ \cite{JJJ}. Furthermore, we establish that if a \vosa is selfdual, regular and of CFT type, then $V$ is $g$ is 
  rational for any automorphism $g$ of $V$ of finite order.

Fix $n\in \frac{1}{T'}\Z_+$ and write $n=\ell+\frac{i}{T'}$ for unique nonnegative integers $\ell$ and $i$ with $0\leq i\leq T'-1$. For integers $0\leq r\leq T'-1$ define
\begin{align*}
\delta_i(r)
&=
\left\{
\begin{array}{lll}
1 & \mbox{if} & r\leq i, \\
0 &\mbox{if} & i<r\leq T'-1.
\end{array}
\right.
\end{align*}
Also set $\delta_i(T')=0$. For $v\in V$ and homogeneous $u\in \VRAST$ define
\begin{align*}
u\circ_{g,n} v 
& = 
\res_{z}Y(u,z)v\frac{(1+z)^{\wt u-1+\delta_i(r)+\ell +\frac{r}{T'}}}{z^{2\ell+\delta_i(r)+\delta_i(T'-r)+1}}
\end{align*}
and
\begin{align*}
u\ast_{g,n}v &=
\left\{
\begin{array}{llll}
\displaystyle{\sum_{m=0}^{\ell}(-1)^m{m+\ell\choose \ell}\res_{z}Y(u,z)v\frac{(1+z)^{\wt{u}+\ell}}{z^{\ell+m+1}}} & \mbox{if} & r=0,\\
\\
0 &\mbox{if} &r>0.
\end{array}
\right.
\end{align*}
Extend both $\circ_{g,n}$ and $\ast_{g,n}$ linearly to obtain bilinear products on all of $V$. Let $\OGNV$ be the linear span of elements  of the form $u\circ_{g,n} v$ and $L(-1)u+L(0)u$ for $u,v\in V$.  Then $V^{r*}$ is a subspace of $\OGNV$ if $r\ne 0.$  Define the linear space $\AGNV$ to be the quotient 
\begin{align*}
\AGNV &= V/ \OGNV.
\end{align*}
Then $\AGNV$ is a quotient of $V^{0*}.$

\begin{thm}\label{Petersen}\cite{P}
	Let $V$  and $g$ be as before, and $M$ be an admissible $g$-twisted $V$-module.
	
	(1)	The product $*_{g,n}$ induces the structure of an associative algebra on $\AGNV$ with identity element $\1+\OGNV$ and central element $\omega+\OGNV$.
	
	(2) The identity map on $V$ induces an onto algebra homomorphism from $\AGNV$ to $A_{g,n-\frac{1}{T'}}(V)$.
	
	(3) The map $u+\OGNV\mapsto o(u)$ gives a representation of the associative algebra $\AGNV$ on $M(m)$
	for $0\leq m\leq n$ where $o(u)=u_{\wt u-1}$ if $u$ is homogeeous.
	
	(4) $M$ is irreducible if and only if $M(n)$ is an irreducible $\AGNV$-module for all $n\geq 0.$
	
	(5) Two admissible $g$-twisted $V$-modules are isomorphic if and only if $M(n)$ and $N(n)$ are isomorphic 
	$\AGNV$-modules for all $n\geq 0.$ 
	
	(6) $V$ is $g$-rational if and only if $\AGNV$ is a finite dimensional semisimple associative algebra for all $n\geq 0.$
\end{thm}

Recall from \cite{DLM6} that if $V$ is a vertex operator algebra with an automorphism $g$ of order $T.$ Then 
for any $n=l+\frac{r}{T}$ with $0\leq r<T,$  $A_{g,n}(V)$ is a quotient algebra of $A_{l}(V)$
where $A_l(V)$ is the associative algebra studied in \cite{DLM5}.  Here is an analogue result for
  \vosa $V.$
  
\begin{lem}\label{relation}
For any integer $n=l+\frac{r}{T'},$  $A_{g,n}(V)$ is quotient algebra of  $A_{\sigma, [n]}(V^{\<\sigma g\>})$ where $\<\sigma g\>$ is the cyclic group generated  by $\sigma g.$
\end{lem} 

\pf  From the definition we know that $O_{\sigma, l}(V^{\<\sigma g\>})$ consists of vectors
$$u\circ_{\sigma,l} v = \res_{z}Y(u,z)v\frac{(1+z)^{\wt u+\ell}}{z^{2\ell+2}},\  (L(0)+L(-1))u$$
for $u,v\in V^{\<\sigma g\>}.$ From the defintion of $A_{g,n}(V)$ we see that 
$u\circ_{g,n} v =u\circ_{\sigma,l} v $
for the same $u,v.$ Since $V^{0*}=V^{\<\sigma g\>}$ and  $*_{g,n}=*_{\sigma, l},$   the result follows. 
\qed

Let $G$ be a finite automorphism group of $V$ such that $\sigma\in G.$
Set $\overline{G}=G/<\sigma>$ and $\bar{g}=g+ <\sigma>\in \overline{G}$ for $g\in G.$ Then $\overline{G}$ is an automorphism 
group of vertex operator algebra $V_{\bar 0}$ isomorphic go $G|_{V_{\bar 0}}$ by identifying $\overline{G}$ with $g|_{V_{\bar 0}}.$ Observe that $V^G=V_{\bar 0}^{\overline{G}}.$ For any subgroup $H$ of $G$ we denote by $\overline{Hu}$ the image of $H$ in $\overline{G}.$ The following lemma is clear from the definition.
\begin{lem}\label{l2} If $M$ is an admissible $g$ or $\sigma g$-twisted $V$-module, then $M$ is a $\bar{g}$-twisted $V_{\bar{0}}$-module.
\end{lem}

The next result seems new and is necessary for applying modular invariance results for \vosa from \cite{DZ1}.
\begin{prop} \label{newrational}
 If a \vosa is self-dual, regular and of CFT type, then $V$ is $g$ is 
rational for any automorphism $g$ of $V$ of finite order.
\end{prop}

\pf This result was obtained previously in  \cite{ADJR} when $V$ is a vertex operator algebra. The proof here is essentialy the same as that given in   \cite{ADJR} with the help from \cite{P}  on the $A_{g,n}(V)$  theory .

Let $G$ be the subgroup of $\Aut(V)$ generated by $g$ and $\sigma.$ Then $G$ is an finite abelian group. As pointed out in \cite{DNR1} that   $V$ is regular if and only if  $V_{\bar 0}$ is regular by using the arguments from \cite{M} and \cite{CM} . Since $\overline G$ is a finite abelian group, $V^G=V^{\overline G}_{\bar 0}$ is also simple, regular and of CFT type. Note that $V_{\bar 0}^{\<g\>}=V^{\<\sigma g\>}_{\bar 0}=V^G.$ As a result,
$V^{\<\sigma g \>}=V_{\bar  0}^{\<\sigma g\>}\oplus V_{\bar 1}^{\<\sigma g\>}$ is a simple, regular and of CFT type. Moreover,
$V^{\<\sigma g\>}$ is $\sigma$-rational.  

From Theorem \ref{Petersen}, $V$ is $g$-rational if and only if $A_{g,n}(V)$
is finite dimensional semisimple associative algebras for all  $n.$ By Lemma \ref{relation}, $A_{g,n}(V)$ is a quotient 
of semisimple associative algebra $A_{l,\sigma}(V^{\<\sigma g\>})$
where  $n=l+\frac{r}{T'},$ and is semisimple for all $n.$  Thus $V$ is $g$-rational.
\qed

To proceed further, we need the action of $\Aut(V)$ on twisted modules \cite{DLM7,DZ1}. 
Let $g, h$ be two automorphisms of $V$ with $g$ of finite order. 
If $(M, Y_M)$ is a weak $\sigma g$-twisted $V$-module, 
there is a weak $h^{-1}\sigma gh$-twisted  $V$-module $(M\circ h, Y_{M\circ h})$ where $M\circ h= M$ as vector spaces and
\begin{equation*}
	Y_{M\circ h}(v,z)=Y_M(hv,z)
\end{equation*}
for $v\in V.$
This defines a right action of $\Aut(V)$ on weak twisted $V$-modules and on isomorphism
classes of weak twisted $V$-modules. Symbolically, we write
\begin{equation*}
	(M,Y_M)\circ h=(M\circ h,Y_{M\circ h})= M\circ h,
\end{equation*}
where we sometimes abuse notation slightly by identifying $(M, Y_M)$ with the isomorphism
class that it defines. Note that if $M$ is an admissible (ordinary, resp.), then $M\circ h$ is admissible (ordinary, resp.).
$M$ is called $h$-stable if $M\circ h$ and $M$ are isomorphic. 

We now define  admissible, ordinary twisted {\em super} $V$-modules.  
An admissible twisted $V$-module $M$ is called an admissible twisted  super $V$-module if $M$ and $M\circ \sigma$ are isomorphic.
In particular, any irreducible $V$-module  is a super module.

\begin{lem}\label{l2.2} If $M$ is an irreducible admissible $g$-twisted $V$-module, 
then $M$ is an irreducible $\bar{g}$-twisted $V_{\bar 0}$-module if $M$ is not $\sigma$-stable, and $M$ is a direct sum of
 two inequivalent irreducible $\bar{g}$-twisted $V_{\bar 0}$-modules if $M$ is $\sigma$-stable.
\end{lem}

\pf
From \cite[Lemma 6.1]{DZ1}, $M$ is $\sigma g$-stable $V$-module. 
If $M\circ\sigma$ and $M$ are not isomorphic, then $M$ is not $g$-stable.
It follows from the proof of  \cite[Theorem 6.1]{DM1} that
$M\circ \sigma$ and $M$ are isomorphic irreducible $\overline{g}$-twisted $V_{\bar 0}$-modules.

If $M$ is $\sigma$-stable,
then there is a linear map $\sigma: M\to M$ such that 
$\sigma Y_M(v,z)\sigma^{-1}=Y_M(\sigma v,z)$ for $v$ in $ V.$
We can  choose $\sigma$ such that $\sigma^2=1.$ Denote the eigenspace with eigenvalue $(-1)^i$ by $M_{\bar i}.$ 
Then $M_{\bar i}$ is irreducible $\overline{g}$-twisted $V_{\bar 0}$-module.  Note that the action of $\sigma$ is not canonical. Replacing $\sigma$ by $-\sigma,$ the new $M_{\bar i}$ will be the old $M_{\overline{i+1}}$ for $i=0,1.$
\qed

\begin{rem}\label{r2.11}  Let $M$ be a nonzero  admissible  $V$-module.  Then $M=M_{\bar 0}\oplus M_{\bar 1}$ is $\sigma$-stable where $M_{\bar i}=\oplus_{m\in \frac{i}{2}\Z_+}M(m)\ne 0$, $i=0,1$ and $\sigma|=_{M_{\bar i}}=(-1)^i.$  Moreover, if $M$ is irreducible, then $M_{\bar 0}$ and $M_{\bar 1}$ are inequivalent irreducible $V_{\bar {0}}$-modules \cite{DNR1}.
\end{rem}

\begin{lem}\label{stable} 
If $M$ be an admissible $g$-twisted $V$-module, then $M\oplus M\circ \sigma$ is $g$-stable super $V$-module.
Moreover, $M$ is $g$-stable if and only if  $M$ is $\sigma$-stable. That is, any admissible super $g$-twisted $V$-module is $g$-stable.
\end{lem}

\pf  From the proof of   \cite[Lemma 2.4]{DNR1},  it is clearly that $M\oplus M\circ \sigma$ is an admissible $g$-twisted super module.   From the proof of \cite[Lemma 6.1]{DZ1},   $M$ is $\sigma g$-stable.  Thus  $M\oplus M\circ \sigma$
is $g$-stable.    It is evident now  that $M$  is  $g$-stable if and only if $M$ is $\sigma$-stable. \qed

\section{Modular Invariance}

We present  the modular invariance property of the trace functions in orbifold theory  for \vosa
from \cite{H, DZ1}. Also see \cite{DLM7, Z}.  We assume in this section  that  $V$ is a simple, regular \vosa of CFT type and $G$ is  any finite automorphism group of $V.$

If $g,h\in G$ commute, $h$ clearly acts on the $g$-twisted modules.
Denote by $\mathscr{M}(g)$ the equivalence classes of irreducible $\sigma g$-twisted super $V$-modules 
and set 
$$\mathscr{M}(g, h)=\{M \in \mathscr{M}(g)\mid M\circ\sigma h \cong M\}.$$ 
Note from Proposition  \ref{newrational} that
 $V$ is $\sigma g$-rational, both $\mathscr{M}(g)$ and $\mathscr{M}(g,h)$ are finite sets.
 For any $M\in \mathscr{M}( g,h),$ there is a $\sigma g$-twisted $V$-module isomorphism
\begin{equation*}
\phi(\sigma h) : M\to M\circ \sigma h
\end{equation*}
such that $\phi(\sigma h)Y_M(v,z)\phi(\sigma h)^{-1}=Y_M(\sigma hv,z)$ for all $v\in V$.
The linear map $\phi(\sigma h)$ is unique up to a nonzero scalar. If $h=\sigma$ we simply take $\phi(1)=1.$
Let  $o(g)=T$.  Then  a $\sigma g$-twisted $V$-module $M$ can decomposes into $M=\oplus_{n\in\frac{1}{T}\Z_+} M_{\l+n}$ (see Theorem \ref{grational}).
For homogeneous  $v\in V$ we set
\begin{equation*}
Z_M(v, (g,h),q)=\tr_{_M}o(v)\phi(\sigma h) q^{L(0)-c/24}=q^{\lambda-c/24}\sum_{n\in\frac{1}{T}\Z_+}\tr_{_{M_{\l+n}}}o(v)\phi(\sigma h)q^{n}.
\end{equation*}
This is a formal power series in variable $q.$
It is proved in  \cite{DZ1, Z, DLM7} that
$Z_M(v,(g,h), q)$ converges to a holomorphic function on the domain $|q|<1.$ We also use $Z_M(v,(g,h),\tau)$ to denote the corresponding holomorphic function. Here and below $\tau$ is in the complex upper half-plane $\H$ and $q=e^{2\pi i\tau}.$
Note that $Z_M(v, (g,h),\tau)$ is defined up to a nonzero scalar. 
Set 
$$Z_M(v,q)=Z_M(v,(g,\sigma),q)=q^{\lambda-c/24}\sum_{n\in\frac{1}{T}\Z_+}\tr_{_{M_{\l+n}}}o(v)q^{n}.$$
  If $v=\1$ is the vacuum then $\ch_M(q)=Z_M(\1,q)$ is  called the formal character or the $q$ character  of $M.$  We aslo denote  $Z_M(1,\tau)$ by  $\chi_M(\tau).$

We recall the  vertex operator superalgebra $(V, Y[~], \1, \tilde{\omega})$ defined in \cite{Z, DZ1}.  
Here $\tilde{\omega}=\omega-c/24$ and
$$Y[v,z]=Y(v,e^z-1)e^{z\cdot \wt v}=\sum_{n\in \Z}v[n]z^{-n-1}$$
for homogeneous $v.$ 
Write
$$Y[\tilde{\omega},z]=\sum_{n\in \Z}L[n]z^{-n-2}.$$
If $v\in V$ is homogeneous with respect to $L[0],$ we denote its weight by $\wt [v].$

Let $P(G)$ denote the ordered commuting pairs in $G.$
For $(g,h)\in P(G)$ and $M\in \mathscr{M}(g,h),$ 
let  $W$ be the vector space spanned by $Z_M(v,(g,h),\tau)$, where each $Z_M$ is garded as a function on $V\times \H.$ 
 Define an action of the modular group $\Gamma=SL(2,\Z)$ on $W$ such that
\begin{equation*}
Z_M|_\gamma(v,(g,h),\tau)=(c\tau+d)^{-{\rm wt}[v]}Z_M(v,(g,h),\frac{a \tau+b}{c\tau +d}),
\end{equation*}
where $ \gamma=\left(\begin{array}{cc}a & b\\ c & d\end{array}\right)\in\Gamma.$
Also, $\Gamma$ acts on the right of $P(G)$ such that $(g, h)\gamma = (g^ah^c, g^bh^d ).$

The following result in \cite{DZ1} requires assumption that $V$ is $g$-rational for all $g\in G.$ This assmption holds
by Proposition  \ref{newrational}.
\begin{thm}\label{minvariance} Let $V$ be a simple, regular, self-dual vertex operator superalgebra of CFT type 
and $G$ a finite automorphism group of $V.$  Then 

(1) There is a representation $\rho: \Gamma\to GL(W)$
such that for $(g, h)\in P(G),$   $\gamma =\left(\begin{array}{cc}a & b\\ c & d\end{array}\right)\in \Gamma,$
and $M\in \mathscr{M}(g,h),$
$$
	Z_{M}|_{\gamma}(v,(g,h),\tau)=\sum_{N\in \mathscr{M}( g^ah^c, g^bh^d)} \gamma_{M,N}^{(g,h)} Z_{N}(v,(g^ah^c,g^bh^d), ~\tau)$$
where $\rho(\gamma)=(\gamma_{M,N}^{(g,h)}).$
That is,
$$Z_{M}(v,(g,h),\gamma\tau)=(c\tau+d)^{{\rm wt}[v]}\sum_{N\in \mathscr{M}( g^ah^c, g^bh^d)}\gamma_{M,N}^{(g,h)} Z_{N}(v,(g^ah^c,g^bh^d), ~\tau).$$

(2) We have $|\mathscr{M}(g,h)|=|\mathscr{M}( g^ah^c, g^bh^d)|$ for any $(g,h)\in P(G)$ and $\gamma\in\Gamma.$

(3) The number of inequivalent irreducible $g$-twisted super $V$-modules is exactly the number of irreducible $g$-stable $V$-modules.  
\end{thm}

 The vector space ${\cal C}(g,h)$ spanned by $Z_M(v,(g,h),\tau)$ for $M\in \mathscr{M}(g,h)$ is called $(g,h)$ confmormal block which has a basis  $Z_M(v,(g,h),\tau)$ for $M\in \mathscr{M}(g,h)$ \cite{DZ1}.   The $\gamma$ gives a linear isomorphsim from  ${\cal C}(g,h) $  to 
${\cal C}(g^ah^c,g^bh^d)$  and (2) follows. Now consider commuting pair $(g\sigma, \sigma)$ and 
matrix $ s=\left(\begin{array}{cc}0 & 1\\ -1 & 0\end{array}\right)\in \Gamma.$  Then $ |\mathscr{M}(g\sigma,\sigma)|
=|\mathscr{M}(\sigma,g\sigma)|.$ Note that $ |\mathscr{M}(g\sigma,\sigma)|$ is exactly the number of inequivalent irreducible $g$-twisted super $V$-modules  and $|\mathscr{M}(\sigma,g\sigma)|$  is exactly the number of  irreducible $g$-stable $V$-modules.

The modular group $\Gamma$ is generated by $S=\left(\begin{array}{cc}0 & -1\\ 1 & 0\end{array}\right)$ and $T=\left(\begin{array}{cc}1 & 1\\ 0 & 1\end{array}\right)$ and the representation $\rho$ is uniquely determined by $\rho(S)$ and  $\rho(T).$  The matrix $\rho(S)$
is called the $S$-matrix of the orbifold theory.

In the rest of this  paper, we assume that $V=\oplus_{n\geq 0}V_n$ is a simple, selfdual  vertex operator superalgebra of CFT type such that
\begin{enumerate}
	\item[(V1)] {\em $G$ is a finite automorphism group of $V$ with  $\sigma\in G$ and $V^G$ is a regular vertex operator algebra,}
	\item[(V2)] {\em The conformal weight of any  irreducible $g$ twisted $V$-module for $g\in G$ except
		$V$ itself is positive.}
\end{enumerate}

Clearly, $V_{\bar 0}$ is an extension of $V^G.$ The following result is immediate from \cite{HKL} and \cite{ABD}.
\begin{thm}\label{M}  
 The $V_{\bar0}$ is regular, selfdual 	vertex operator algebra of CFT type.
\end{thm} 

We point out  that if $G$ is solvable, $V$ is regular if and  only if $V^G$  is regular  from \cite{ABD, CM, HKL, M}.

\section{Classification of irreducible modules for $V^G$}
In this section we classify the irreducible $V^G$ modules  and show that any irreducible $V^G$ module occurs in an irreducible $g$-twisted $V$-module for some $g\in G.$  In the case $V$ is a vertex operator algebra, this result was obtained previously in \cite{DRX1}. We also establish that the orbifold trace functions are modular forms on a congruence 
subgroup. generalizing a similar result when $V$ is a vertex operator algebra  \cite{DR}.

Let $M=(M,Y_M)$ be an irreducible $g$-twisted $V$-module.
We define a subgroup $G_M$ of $G$ consisting of $h\in G$ such that  $M$ is $h$-stable.  Note that $G_M$ is a subgroup of $C_G(g).$  We claim that $\sigma g$ lies in $G_M.$ 
From Theorem \ref{grational}, $M$ has a decomposition $M=\oplus_{n\in \frac{1}{T'}\Z_+}M(n)$ and $M(0)\ne 0$, where $o(\sigma g)=T'.$
We define $\phi(\sigma g)$ acting on $M_{\lambda+\frac{n}{T'}}$ as $e^{2\pi in/T'}$ for all $n.$ It is easy to check that
$$
\phi(\sigma g)Y_M(v,z)\phi(\sigma g)^{-1}=Y_M(\sigma gv,z)
$$
for  $v\in V.$ This implies $\sigma g \in G_M.$ 

The $M$ is a projective
$G_M$-modue such that  $h\in G_M$ acts on $M$ as  $\phi(h)$ satisifying
$$
\phi(h)Y_M(v,z)\phi(h)^{-1}=Y_M(hv,z)
$$
for $v\in V.$
Let $\a_M$  be the corresponding 2-cocycle in $C^2(G_M,\C^{\times}).$ Then $\phi(h)\phi(k)=\a_M(h,k)\phi(hk)$
for all $h,k\in G_M.$ 
We may  assume that $\alpha_M$ is unitary  in the sense that there is a fixed positive integer $n$ such that  $\alpha_M(h,k)^n=1$ for all $h,k\in G_M$ \cite{C}.  
Let $\C^{\a_M}[G_M]=\oplus_{h\in G_M} \C\hat{ h}$ be the twisted group algebra with product $
\hat h\hat k=\alpha_M(h,k)\hat{hk}.$ 
Then $\C^{\a_M}[G_M]$ is a semisimple associative algebra and $M$ is a $\C^{\a_M}[G_M]$-module such that $\hat h$ acts as $\phi(h).$
We will denote the irreducible characters of $\C^{\a_M}[G_M]$ by $\Lambda_M.$ In the case $\sigma\in G_M,$ 
 $\Lambda_M= \Lambda_{M}^0\cup  \Lambda_M^1$ where $\Lambda_M^i$ consists of the irreducible characters 
$\lambda$ such that $\sigma=(-1)^i$ on the corresponding irreducible $\C^{\a_M}[G_M]$-module as $\sigma=\hat{\sigma}$ is a central element of order 2.
For $\lambda\in  \Lambda_M$ we denote the corresponding irreducible  $\C^{\a_M}[G_M]$-module by $W_{\lambda}.$

Since the actions of $\C^{\alpha_M}[G_M]$ and $V^G$ commute on $M,$ we  now study a Schur-Weyl duality type decomposition of $M.$  Note that  $M$ is a semisimple $\C^{{\a}_M}[G_M]$-module as $\C^{{\a}_M}[G_M]$
is a finite dimensional semisimple associative algebra. For $\lambda\in \Lambda_M,$  let $M^{\lambda}$ be the sum of simple $\C^{{\a}_M}[G_M]$-submodules of $M$ isomorphic to $W_{\l}.$ Then
$$M=\oplus_{\lambda\in \Lambda_M}M^{\lambda}=\oplus_{\lambda\in \Lambda_M}W_{\lambda}\otimes M_{\l}$$
where $M_{\l}=\Hom_{\C^{{\a}_M}[G_M]}(W_{\l},M)$ is the multiplicity
of $W_{\l}$ in $M.$ Then $M_{\lambda}$ is a $V^{G_M}$-module such that for $u\in V^{G_M},$ $m\in\Z$ and
$f\in M_{\lambda},$ $(v_mf)(w)=v_mf(w).$ 
As in \cite{DLM2, DY, DRX1}, we can,
in fact, realize $M_{\l}$ as a subspace of $M$ in the following way:
Let $w\in W_{\l}$ be a fixed nonzero vector. 
Then we can identify
$\Hom_{\C^{{\a}_M}[G_M]}(W_{\l},M)$ with the subspace
$$\{f(w) |f\in \Hom_{\C^{{\a}_M}[G_M]}(W_{\l},M)\}$$
of $M^{\l}.$ It is easy to prove directly that this subspace  is a $V^{G_M}$-module.

Now we consider vertex operator algebra $V_{\bar 0}$ and the group $\overline{G}.$ By Lemma \ref{l2.2}, $M$ is an 
irreducible $\overline{g}$-twisted $V_{\bar 0}$ module if  $M$ is not $\sigma$-stable and 
$M=M_{\bar 0}\oplus M_{\bar 1}$ is a direct sum of two inequivalent irreducible  $\overline{g}$-twisted $V_{\bar 0}$ modules.
For any  $\overline{g}$-twisted $V_{\bar 0}$ module $W$ we let $\overline{G}_W$ be the subgroup of $\overline{G}$ consisting of
$\bar k\in\overline{G}$ such that $W\circ \bar k\cong W.$ 

\begin{lem}\label{mainlemma} Let $M$ be an  irreducible $g$-twisted $V$ module and $N$ irreducible $h$-twisted $V.$ Then

(1) If  $M,N$ are not $\sigma$-stable, then $M,N$ are isomorphic $\overline{g}$-twisted $V_{\bar 0}$ if and only if
$g=h$ and either $M\cong N$ or $M\circ \sigma\cong N$ as $g$-twisted $V$-modules.

(2) If $M$ is $\sigma$-stabe and $N$ is not $\sigma$-stable then $M_{\bar i}$ and $N$ are not isomorphic 
$\overline{g}$-twisted $V$-modules.

(3)  If  $M,N$ are $\sigma$-stable, then $M_{\bar i}$ and $N_{\bar j}$ are  isomorphic $\overline{g}$-twisted $V_{\bar 0}$-modues  
for some $i,j$ if and ony if $M\cong N.$

(4) If $M$ is not $\sigma$-stable then $\overline{G}_M=\overline{G_M}$ is isomorphic to $G_M.$ 

(5) If $M$ is $\sigma$-stable, then $\overline{G}_{M_{\bar 0}}=\overline{G}_{M_{\bar 1}}$ is a subgroup of $\overline{G_M}.$
Moreover, 
the index $[\overline{G_M}: \overline{G}_{M_{\bar 0}}]$ is 1 or 2, and the index equals to 1  if and only if $\phi(k)\sigma=
\sigma\phi(k)$ for any $k\in G_M.$ 
\end{lem}
\pf  Some category theory is necessay to prove the result. Recall some basics on  category theory from \cite{KO, EGNO}.  An object $A$ in a braided fusion category $\CC$ is called regular commutative algebra  if there are morphisms
$\mu: A\boxtimes A\to A$ and $\eta: {\bf 1_\CC}\to A$ such that $\mu\circ(\mu\boxtimes {\rm id}_A)\circ \alpha_{A,A,A} = \mu\circ({\rm id}_A\boxtimes\mu)$,  $\mu\circ(\eta\boxtimes {\rm id}_A)={\rm id}_A =\mu\circ({\rm id}_A\boxtimes\eta),$ $\mu = \mu\circ c_{A,A},$ 
 and
$\dim\hom({\1}_\CC,A)=1$ where $\alpha_{A,A,A}: A\boxtimes (A\boxtimes A)\to (A\boxtimes A)\boxtimes A$ is the associative isomorphism,  and $c_{A,A}: A\boxtimes A\to A\boxtimes A$ is the braiding isomorphism. A left $A$-module $N$ is an object in $\CC$ with a morphism $\mu_N: A\boxtimes N\to N$ such that $\mu_N\circ (\mu\boxtimes \id_N)\circ \alpha_{A,A,N}=\mu_N\circ (\id_A\boxtimes \mu_N).$ We denote the left $A$-module category by $\CC_A.$ Then  $\CC_A$ is a fusion category with tensor product $N_1\boxtimes_A N_2$. An $A$-module $N$ is called local if
$\mu_N\circ c_{N,A}\circ c_{A,N}=\mu_N.$  The local $A$-module categoy  $\CC_A^0$ is a braided fusion category. Moreover, if $\CC$ is modular tensor category, so is $\CC_A^0$ \cite{KO}. For 
any $N\in \CC_A$ and $X\in\CC,$ $N\boxtimes X\in \CC_A.$ From \cite{Hu},  the $V^G$-module category ${\cal C}_{V^G}$  is a modular tensor category with tensor product $\boxtimes$ as $V^G$ is regular.  Moreover, $V_{\bar 0}$ is a regular commutative algebra in $\CC_{V^G}$ \cite{HKL}. From the discussion above, $(\CC_{V^G})_{V_{\bar 0}}$  is a fusion category  and  $(\CC_{V^G})_{V_{\bar0 }}^0$  is exactly $\CC_{V_{\bar 0}}.$ We also know from \cite{DLXY} that every simple object in  $(\CC_{V^G})_{V_{\bar 0}}$  is isomorphic to an irreducible $\bar h$-twisted $V_{\bar 0}$-module, and the tensor products of an twisted $V_{\bar 0}$-module and a twisted $V_{\bar 0}$-module in the vertex operator algebra setting \cite{X} and the categorical setting \cite{KO} are the same. 

It is well known that $V_{\bar 1}$ is a simple current  in $\CC_{V_{\bar 0}}=(\CC_{V^G})_{V_{\bar0 }}^0$ (cf. \cite{DNR1}). Then  $V_{\bar 1}$  is also a simple current in category $(\CC_{V^G})_{V_{\bar0 }}$ by using the associativity of the fusion product $\boxtimes_{V_{\bar 0}}.$  So for any irreducilbe $\overline{g}$-twisted $V_{\bar 0}$-module $W,$ $V_{\bar 1}\boxtimes_{V_{\bar 0}} W$ is again an irreducible  $\overline{g}$-twisted $V_{\bar 0}$-module for any $g\in G.$ We will use this fact in the following proof.

(1) It is good enough to prove the ``if'' part.  Since $V_{\bar 1}$ is a simple current,  $V_{\bar 1}\boxtimes_{V_{\bar 0}} M\cong M $ is isomorhic $V_{\bar 1}\boxtimes_{V_{\bar 0}} N\cong N$ as $\overline{g}$-twisted $V_{\bar 0}$-modules.
Without loss we can assume that $M=N$ as $\overline{g}$-twisted $V_{\bar 0}$-modules.
Let $M=(M,Y_1)$ and $N=(M,Y_2)$ with $Y_1(v,z)=Y_2(v,z)$ for $v\in V_{\bar 0}.$
Note that the space of intertwining operators of type   $\left(\begin{array}{c}
M\\
V_{\bar 1}\  M
\end{array}\right)$   is one dimensional and spanned by $Y_1,$ and $Y_2$ is also intertwining operator
of the same type. Thus there exist a nonzero constant $\lambda$ such that  $Y_2(v,z)=\l Y_1(v,z)$ for all $v\in V_{\bar 1}.$ Let $u\in V_{\bar 1}^r$ and $v\in  V_{\bar 1}.$ Using the associativity for vertex operator $Y_2(u,z_1)$ and $Y_2(v,z_2)$ for $u,v\in V_{\bar 1}$ we see that $\lambda=\pm 1.$ If $\lambda=1$ then $N=M$, and if $\lambda=-1$
then $N=M\circ \sigma$ as $g$-twisted $V$-module. In particular, $h=g.$ 

(2) If $M_{\bar 0}$ and $N$ are isomorphic $\overline{g}$-twisted $V_{\bar 0}$-modules, then 
$V_{\bar 1}\boxtimes_{V_{\bar 0}} M_{\bar 0}\cong M_{\bar 1}$ and $V_{\bar 1}\boxtimes_{V_{\bar 0}}  N\cong N$ are isomorphic $\overline{g}$-twisted $V_{\bar 0}$-modules. As a result, $M_{\bar 0}$ and $M_{\bar 1}$ are isomorphic $\overline{g}$-twisted $V_{\bar 0}$-modules. This contradicts to Lemma \ref{l2.2}.

(3) Again we only need to prove the ``if'' part. Without loss we can assume that $M_{\bar 0}$ and $N_{\bar 0}$ are 
 $\overline{g}$-twisted $V_{\bar 0}$-modules as the actions of $\sigma$  on $M, N$ can be replaced by $-\sigma$ if necessay. In this case, $V_{\bar 1}\boxtimes_{V_{\bar 0}}  M_{\bar 0}\cong M_{\bar 1}$ and $V_{\bar 1}\boxtimes_{V_{\bar 0}}  N_{\bar 0}\cong N_{\bar 1}.$ This impies that $M,N$ are isomorphic $\overline{g}$-twisted $V_{\bar 0}$-modules. As in (1) we can assume that $M=N$
 as $V_{\bar 0}$-modules with module maps $Y_1,Y_2$ such that (i) $Y_1(v,z)=Y_2(v,z)$ for $v\in V_{\bar 0},$
 (ii) There exists nonzero constant $\lambda$ such that $Y_2(u,z)w=\lambda Y_1(u,z)w$ for $u\in V_{\bar 1}$ and
 $w\in M_{\bar 0},$ (iii) $Y_2(u,z)w=\lambda^{-1} Y_1(u,z)w$ for $u\in V_{\bar 1}$ and $w\in M_{\bar 1}.$ Define
 $\phi:  M\to M$ such that $\phi(w^0+w^1)=w^0+\lambda w^1$ for $w^i \in M_{\bar i}.$ We claim that $\phi Y_1(u,z)w=
 Y_2(u,z)\phi(w)$ for $u\in V$ and $w\in M.$ This is cleay if $u\in V_{\bar 0}.$  Now let $u\in V_{\bar 1}$ and
 $w=w^0+w^1$ wth $w^i\in M_{\bar i}.$ Then 
 $$\phi Y_1(u,z)w=\lambda Y_1(u,z)w^0+Y_1(u,z)w^1=Y_2(u,z)\phi(w^0)+ Y_2(u,z)\phi(w^1)=Y_2(u,z)\phi(w).$$
 That is, $\phi$ gives an isomorphism between $g$-twisted $V$-modules $(M,Y_1)$ and $(M,Y_2).$ 

(4) Clearly, $\overline{G_M}$  is a subgroup of $\overline{G}_M.$ Let $k\in G$ such that $\bar k\in \overline{G}_M.$ Then 
$M\circ \bar k$ and $M$ are isomorphic  $\overline{g}$-twisted $V_{\bar 0}$-modules. Equivalently, $M\circ k$ and $M$ are
isomorphic  $\overline{g}$-twisted $V_{\bar 0}$-modules. It follows from (1) that $M\circ k$ is a $g$-twisted $V$-module
 isomorphic to either $M$ or  $M\circ \sigma.$ If  $M\circ k\cong M$ then $k\in G_M$ and  $\bar k\in \overline{G_M}.$ If $M\circ k\cong M\circ \sigma$ then $\sigma k\in G_M$ and $\bar k=\overline{\sigma k}\in \overline{G_M}.$ Thus 
 $\overline{G_M}=\overline{G}_M.$ Since $\sigma$ is not an element of $G_M,$ the map $h\mapsto \bar h$ gives a group isomorphism from $G_M$ to
 $\overline{G_M}.$

(5) We first  prove that that $\overline{G}_{M_{\bar i}}$ is contained in $\overline{G_M}.$  Without loss we can take $i=0.$ Let $k\in G$ such that $M_{\bar 0}\circ \bar k\cong M_{\bar 0}.$ From the proof of (3) we see that 
$M\circ k$ and $M$ are isomorphic $g$-twisted $V$-modues. Thus $k\in G_M$ and $\bar k\in \overline{G_M}.$
It is clear that $M_{\bar 0}\circ \bar k\cong M_{\bar 0}$ if and only if $M_{\bar 1}\circ \bar k\cong M_{\bar 1}.$ This implies that $\overline{G}_{M_{\bar 0}}=\overline{G}_{M_{\bar 1}}.$ 

If $\overline{G}_{M_{\bar 0}}$ is a proper subgroup of $\overline{G_M}$ we claim that $[\overline{G_M}:\overline{G}_{M_{\bar 0}} ]=2.$ Let $G_{M_{\bar i}}=\{k\in G_M|M_{\bar i}\circ \bar k\cong M_{\bar i}\}.$ Then $G_{M_{\bar i}}$ is a proper subgroup of $G_M$ and $\overline{G_{M_{\bar 0}}}=\overline{G}_{M_{\bar 0}}.$ 
So it is good enough to show that 
$[G_M:G_{M_{\bar 0}}]=2,$ or $G_M\setminus G_{M_{\bar 0}}=G_{M_{\bar 0}}k$  for any fixed $k\in G_M\setminus G_{M_{\bar 0}}.$  Clearly $G_{M_{\bar 0}}k$ is a subset of $G_M\setminus G_{M_{\bar 0}}.$
Now let  $h\in G_M\setminus G_{M_{\bar 0}}.$ Then $M_{\bar i}\circ \bar h\cong M_{\overline{i+1}}$ 
for $i=0.1.$ This  forces   $hk^{-1}\in G_{M_{\bar 0}}$ and  $h\in G_{M_{\bar 0}}k.$ 

Finally, $\overline{G}_{M_{\bar 0}}=\overline{G_M}$ if and only if $\phi(k)M_{\bar i}=M_{\bar i}$ for all $k\in G_M.$ This is equivalent that $\sigma \phi(k)=\phi(k)\sigma$ for all $k\in G_M.$
\qed

\begin{rem}It happens that $G_{\bar M_{\bar 0}}$ is a proper subgroup of $G_M$ when $M$ is $\sigma$-stable. Here is an example.
Let $V=V_{\Z\alpha}$ the holomorphic lattice \vosa associated to unimodular lattice $\Z\alpha$ with $(\alpha,\alpha)=1$ \cite{B,FLM}. Then $V_{\bar 0}=V_{2\Z\alpha}$ and $V_{\bar 1}=V_{2\Z\alpha+\alpha}.$ Let $\theta$ be an automorphism of $V$ such that 
$$\theta(\alpha(-n_1)\cdots \alpha(-n_s)e^{n\alpha}=(-1)^s\alpha(-n_1)\cdots \alpha(-n_s)e^{-n\alpha}$$
where $n_1,...,n_s$ are positive integers and $n\in \Z.$ Let $G$ be an automorphism group of $V$ generated by
$\sigma$ and $\theta.$ $V$ has a unique $\sigma$-twisted  $V$-module $M=V_{\Z\alpha+\alpha/2}$ which is $\sigma$-stable such that  $M_{\bar 0}=V_{2\Z\alpha+\alpha/2}$ and  $M_{\bar 1}=V_{2\Z\alpha-\alpha/2}.$
We know from \cite{DN} that $M_{\bar 0}\circ \theta\cong  M_{\bar 1}.$ 
\end{rem}

This is the frist result towards to a classifiction of irreducible $V^G$-modules.
 \begin{thm}\label{mthm1} Let $M$ be an irreducible $g$-twisted $V$-module. 
	
	(1) $M^{\l}$ is nonzero for any $\l\in \Lambda_{M}.$
	
	(2)  Each  $M_{\l}$ is an irreducible $V^{G_M}$-module. 
	
	(3) $M_{\l}$ and $M_{\gamma}$ are equivalent $V^{G_M}$-modules if and only
	if $\l=\gamma.$
	\end{thm}

\pf   The proof of Theorem is similar to that of \cite[Theorem 5.4]{DY} by using the associative algebras $A_{g,n}(V)$
from Section 3.
\qed

The classification of irreducible $V^G$-modules requires to investigate the connection between the $G$-orbits on
the set of all inequivalent irreducible twisted  $V$-modules and the  $\overline{G}$-orbits on
the set of all inequivalent irreducible twisted  $V_{\bar 0}$-modules which appear in irreducible twisted  $V$-modules.
Let $ {\cal S}=\cup_{g\in G} {\cal S}(g)$ where $ {\cal S}(g)$ is the set of inequivalent of irreducible $g$-twisted $V$-modules. Then $G$ acts on $ {\cal S}$  and ${\cal S}=\cup_{j\in J}M^j\circ G$ is a union of disjoint orbits. It is clear that $M,N$ are isomorphic $V^G$-modules if they are in the same $G$-orbits. So it is enough to consider 
  $M^j$ for the purpose of classification of irreducible $V^G$-modules. 
  
 Let 
 $${\cal S}(g)=\{N^i, N^i\circ \sigma, N^k|i=0,...,q, k=q+1,...,p\}$$
 such that $N^i$ for $i\leq q$ are irreducible  $g$-twisted  $V$-modules which are not $\sigma$-stable and 
  $N^k$ for $k>q$ are irreducible  $g$-twisted  $V$-modules which are $\sigma$-stable.  Since $N^i$ and $ N^i\circ \sigma$ are equivalent irreducible $\bar g$-twisted $V_{\bar 0}$-modules, 
  $$\overline{{\cal S}}(\bar g)=\{N^i,N^k_{\bar s}|i=0,...,q, k=q+1,...,p, s=0,1\}$$
  is the set of inequivalent irreducible $\bar g$-twisted $V_{\bar 0}$-modules coming from irreducible $g$-twisted $V$-modules.  We will prove, in fact, later that $\overline{{\cal S}}(\bar g)$ gives all inequivalent irreducible $\bar g$-twisted $V_{\bar 0}$-modules. Also set  $\overline{{\cal S}}=\cup_{\bar g\in \overline{G}}\overline{{\cal S}}(\bar g).$ 
 Then $\overline{G}$ acts on $\overline{{\cal S}}$ in the same way.  In order to determine the $\overline{G}$-orbits of  $\overline{{\cal S}}$ we write $J=J_1\cup J_2\cup J_3$ where $J_1$ consists of $j$ such that $M^j$ is not $\sigma$-stable, $J_2$  consists of $j$ such that $M^j$ is $\sigma$-stable and $\overline{G}_{M^j_{\bar 0}}=\overline{G_{M^j}},$ 
 and $J_3$  consists of $j$ such that $M^j$ is $\sigma$-stable and $\overline{G}_{M^j_{\bar 0}}\ne \overline{G_{M^j}}$ 
(see Lemma  \ref{mainlemma}).

The following is a key lemma.
\begin{lem}\label{mainlemma1} The set 
$$S=\{M^{j_1}, M^{j_2}_{\bar 0},M^{j_2}_{\bar 1}, M^{j_3}_{\bar 0}|j_1\in J_1, j_2\in J_2, j_3\in J_3 \}$$
gives a complete list of orbit representatives of $\overline{G}$-set  $\overline{{\cal S}}.$ 
  \end{lem}
 \pf Note that if $j\in J_1,$ then $M^j\circ \overline{G}=M^j\circ  G$  by regarding $M^j\circ k$ as an irreducible onject in
$\overline{{\cal S}}.$ If  $j\in J_2$  then  the irreducible objects of $\overline{{\cal S}}$ coming from $M^j\circ  G$ 
are $M^j_{\bar 0} \circ \overline{G}$ and $M^j_{\bar 1} \circ \overline{G}.$ In this case, $M^j_{\bar 0}$ and $M^j_{\bar 1}$ are not in the same $\overline{G}$-orbit by Lemma \ref{mainlemma}. Now we assme that $j\in J_3.$ By Lemma \ref{mainlemma},
 $M^j_{\bar 0}$ and $M^j_{\bar 1}$ are in the same $\overline{G}$-orbit and $[\overline{G}_{M^j}: \overline{G}_{M^j_{\bar 0}}]=2.$
 This imples that $|M^j_{\bar 0}\circ \overline{G}|=2|M^j\circ G|.$ Clearly, $M^j_{\bar 0}\circ k$ is an irreducible
 $\overline{k^{-1}gk}$-twisted $V_{\bar 0}$-module contained in $M^j\circ k\in M^j\circ G$ for $k\in G.$ Since $M^j\circ G$
 gives exactly $2|M^j\circ G|$ irreducible objects in $\overline{{\cal S}},$  we see that all the irreducible objects in $\overline{{\cal S}}$ coming from  $M^j\circ G$ form a single $\overline{G}$-orbit $M^j_{\bar 0}\circ \overline{G}.$ 
  The proof is complete.
 \qed
  
  We are now ready to classify the irreducible $V^G$-modules which appear in irreducible $g$-twisted $V$-modules 
  for $g\in G.$
  
  \begin{thm}\label{mt2}  We have 
  	
  	(1) For $j\in J$ and $\lambda\in \Lambda_{M^j}, $ $ M^j_\l$ is an irreducible $V^G$-module,
   	
  	(2) Let $i,j\in J$ and $\l\in \Lambda_{M^i}, \mu\in \Lambda_{M^j}.$ Then $M^i_\l$ and $M^j_{\mu}$ are isomorphic
  	$V^G$-modules if and only if $i=j, \l=\mu.$
  	\end{thm}
 
 \pf For any $X\in S,$ we have decomposition
 $$X=\oplus_{\bar \l\in \bar \Lambda_X}W_{\bar \l}\otimes X_{\bar \lambda}$$
 where  $\bar \Lambda_X$ is the irreducible characters of the twisted group algebra $\C^{\bar \alpha_X}[\overline{G}_X],$
 $W_{\bar \l}$ is the irreducibe $\C^{\bar \alpha_X}[\overline{G}_X]$-module affording to the character $\bar\l,$ and 
 and $X_{\bar \lambda}$ is the multiplicity space of $W_{\bar \l}$ in $X.$
 By \cite[Theorem 2]{MT} (also see \cite{DLM3, DY}) we know that $\{X_{\bar\l}|X\in S, \bar\l\in \bar \Lambda_X\}$ are inequivalent 
 irreducible $V^G=V^{\overline{G}}_{\bar 0}$-modules. We need to prove $\{X_{\bar\l}|X\in S, \bar\l\in \bar \Lambda_X\}=
 \{M^j_{\l} |j\in J, \l\in G_{M^j}\}.$ That is, we need to identify
 $X_{\bar\lambda}$ with $M^j_\mu$ 
 for $\bar \lambda\in \bar \Lambda_X$ and $\mu \in \Lambda_{M^j}.$  Let $M=M^j$ and $g=g_j.$
 
 1) Let $j\in J_1.$ Then $G_M, \overline{G}_M$ are isomophic and $\alpha_M=\bar \alpha_M$ by Lemma \ref{mainlemma}. Then $\Lambda_M=\bar \Lambda_M$ and the identification is clear. 
 
 2) Let $j\in J_2.$ Then $M=M_{\bar 0}\oplus M_{\bar 1},$ and 
 $M_{\bar 0}$ and $M_{\bar 1}$ are inequivalent irreducible $\overline{g}$-twisted $V_{\bar 0}$-modules where
 $\sigma =(-1)^i$ on $M_{\bar i}$  for $i=0,1$ by Lemma \ref{l2.2}.  In this case, $\sigma$ is a central element of
 $\C^{\alpha_M}[G_M]$ and $\Lambda_M=\Lambda_M^0\cup \Lambda_M^1$ where $\lambda\in \Lambda_M^i$
 if and only if $\sigma|_{W_\l}=(-1)^i.$  One can easily see 
 $$M_{\bar i}=\oplus_{\lambda\in \Lambda_M^i}W_{\l}\otimes M_{\l}.$$
 Let $I_i=\C^{ \alpha_{M}}[G_M](\sigma-(-1)^i)$ be the two-sided  ideal of of $\C^{ \alpha_{M}}[G_M].$
 Then 
 $$\C^{ \alpha_{M}}[G_M]\cong \C^{ \alpha_{M}}[G_M]/I_0\oplus \C^{ \alpha_{M}}[G_M]/I_1.$$
 Note the $\Lambda_M^i$ is exactly the set of irreducible characters of 
 semisimple associative algebra $\C^{ \alpha_{M}}[G_M]/I_i.$
 
 On the other hand, $\overline{G}_{M_{\bar i}}=\overline{G_{M}}$ (see Lemma \ref{mainlemma})  and 
 $$M_{\bar i}=\oplus_{\bar\lambda\in \Lambda_{M_{\bar i}}}W_{\bar \l}\otimes (M_{\bar i})_{\bar\l}.$$
 We claim  that $\C^{\bar \alpha_{M_{\bar i}}}
[\overline{G}_{M_{\bar i}}]$ is isomorphic to $\C^{ \alpha_{M}}[G_M]/I_i.$ This requires a better understanding of algebra $\C^{\bar \alpha_{M_{\bar i}}}
[\overline{G}_{M_{\bar i}}]$ and its action on $M_{\bar i}.$  Let $K=\<\sigma\>$ and
$\{h_a|a\in G_M/K\}$ be a
set of coset representatives of $K$ in $G_M$ such that $h_{K}=1.$ Then $ \overline{G_{M}}=\{\overline{h_a}|a\in G_M/K\}.$ We now define the projective action of $\overline{G}_{M_{\bar i}}$ on $M_{\bar i}$ by $\phi(\overline{h_a})=\phi(h_a).$ Then 
$$\phi(\overline{h_a})\phi(\overline{h_b})=\phi(h_a)\phi(h_b)=\alpha_M(h_a,h_b)\phi(h_ah_b)$$
for $a,b\in G_M/K.$ Note that $h_ah_b=\sigma^r h_{ab}$ for some $r=0,1$ as $\overline{h_a}\,\overline{h_b} =\overline{h_{ab}}.$ So 
$$\phi(h_ah_b)=\alpha_M(\sigma^r,h_{ab})^{-1}\sigma^{r}\phi(h_{ab})=\alpha_M(\sigma^r,h_{ab})^{-1}(-1)^{ir}\phi(h_{ab})$$
where we have used $\sigma$ for $\phi(\sigma).$ This gives
$$\phi(\overline{h_a})\phi(\overline{h_b})=\bar\alpha_{M_{\bar i}}(\overline{h_a},\overline{h_b})\phi(\overline{h_{ab}})$$
with $\bar\alpha_{M_{\bar i}}(\overline{h_a},\overline{h_b})=(-1)^{ri}\alpha_M(h_a,h_b)\alpha_M(\sigma^r,h_{ab})^{-1}.$ 
It is clear now that the action of algebra $\C^{\bar \alpha_{M_{\bar i}}}
[\overline{G}_{M_{\bar i}}]$ on $M_{\bar i}$ is exactly the action of $\C^{ \alpha_{M}}[G_M]$ with $\sigma=(-1)^i.$ 
That is, $\C^{\bar \alpha_{M_{\bar i}}}
[\overline{G}_{M_{\bar i}}]$ is isomorphic to $\C^{ \alpha_{M}}[G_M]/I_i.$ 
As a result, we see that 
$$\{M_{\lambda}|\lambda\in  \Lambda_M^i\}=\{(M_{\bar i})_{\bar\l}| \bar\lambda\in \Lambda_{M_{\bar i}}\}$$
for $i=0,1.$ 

3) Let $j\in J_3.$ Again  $M=M_{\bar 0}\oplus M_{\bar 1}$ but $M_{\bar 0}$ and $M_{\bar 1}$ are equivalent $V^G$-modules.  In this case, $\sigma$ is not a central element of $\C^{\alpha_M}[G_M].$ Recall from the proof of Lemma \ref{mainlemma} that $G_{M_{\bar i}}=\{h\in G_M|M_{\bar i}\circ \bar h\cong M_{\bar i}\},$ $G_{M_{\bar 0}}=
G_{M_{\bar 1}},$  $G_M=G_{M_{\bar 0}}\cup G_{M_{\bar 0}}k$ for any $k\notin G_{M_{\bar 0}}$ and 
$\overline{G}_{M_{\bar 0}}=\overline{ G_{M_{\bar 0}}}.$ Note from 2) that $\C^{\bar \alpha_{M_{\bar i}}}[\overline{G}_{M_{\bar i}}]$ is isomorphic to $\C^{\alpha_M}[G_{M_{\bar 0}}]/I_i$ where $I_i=\C^{\alpha_M}[G_{M_{\bar 0}}](\sigma-(-1)^i)$ for $i=0,1$ 
and $\C^{\alpha_M}[G_{M_{\bar 0}}]\cong \C^{\alpha_M}[G_{M_{\bar 0}}]/I_0\oplus \C^{\alpha_M}[G_{M_{\bar 1}}]/I_1.$ Let $W$ be any irreducible $\C^{\alpha_M}[G_{M_{\bar 0}}]$-module such that $\sigma$ acts as 1. Then the induced module  $W^e=\C^{\alpha_M}[G_M]\otimes_{\C^{\alpha_M}[G_{M_{\bar 0}}]} W=W+\hat kW$ is an irreducible $\C^{\alpha_M}[G_M]$-module and any irredcuble $\C^{\alpha_M}[G_M]$ can be abtained in this way.
We have  decompsotions
 $$M_{\bar 0}=\oplus_{\bar\l\in\Lambda_{M_{\bar 0}}}W_{\bar \lambda}\otimes (M_{\bar 0})_{\bar \l},\ 
 M_{\bar 1}=\phi(k)M_{\bar 0}=\oplus_{\bar \l\in\Lambda_{M_{\bar 0}}}\phi(k)W_{\bar \lambda}\otimes (M_{\bar 0})_{\bar \l}$$
 where $\Lambda_{M_{\bar 0}}$ is the set of irreducible characters of $\C^{\alpha_M}[G_{M_{\bar 0}}]/I_0.$
So 
$$M=\oplus_{\bar \l\in \Lambda_{M_{\bar 0}}}(W_{\bar\l}+\phi(k)W_{\bar \lambda})\otimes (M_{\bar 0})_{\bar \l}=\oplus_{\lambda \in \Lambda_M}W_{\l}\otimes M_{\l}$$
where $W_{\lambda}=W_{\bar\l}+\phi(k)W_{\bar \lambda},$ $\lambda\in \Lambda_M$ is the corresponding character
of $\C^{\alpha_M}[G_M], $ and $(M_{\bar 0})_{\bar \l}=M_\l.$

From 1)-3) we see that 
$$\{M^j_{\l} |j\in J, \l\in G_{M^j}\}=\{X_{\bar \l}|X\in S, \bar\l\in \overline{G}_X\},$$
as desired.
\qed

 This result gives a complete classification of irreducible $V^G$-modules occurring in the irreducible $g$-twisted $V$-modules for $g\in G.$ 
We now use the modular invariance to prove that these irreducible $V^G$-modules are  all the irreducible $V^G$-modules .
\begin{thm}\label{MTH1}  Let $V$ and $G$ be as before. Then $\{M^j_{\l} |j\in J, \l\in G_{M^j}\}$ gives a complete
	list of inequivalent irreducible $V^G$-modules.
\end{thm}

\pf By Theorem \ref{mt2}, it suffices to show that every  irreducible $V^G$-module appears in an 
irreducible $g$-twisted $V$-module $M.$

By Theorems \ref{mthm1}, 
$V=\oplus_{\lambda\in \Lambda_{G}}W_{\lambda}\otimes V_{\lambda},$ where $\Lambda_V=\irr(G)$ is the set of 
irreducible characters of $G$ and $V_{\lambda}$ is an irreducible $V^G$-module.
From the definition of $Z_{V^G}(v,\tau)$ and the orthogonality property of the irreducible characters of $G$, we know that for $v\in V^G$,
$$\frac{1}{|G|}\sum_{g\in G}Z_V(v,(\sigma,\sigma g),\tau)
=\frac{1}{|G|}\sum_{g\in G}\sum_{\lambda\in \Lambda_{V}}(\tr_{W_{\lambda}}g)Z_{V_{\lambda}}(v,\tau)
=Z_{V^G}(v,\tau).$$
Applying Theorem \ref{minvariance} (1) with $\gamma=S=\left(\begin{array}{cc} 0& -1\\ 1 & 0\end{array}\right)$
gives
$$Z_{V^G}(v,-\frac{1}{\tau})=\frac{\tau^{\wt[v]}}{|G|}\sum_{g\in G}\sum_{M\in  \mathscr{M}(\sigma g, \sigma)} S_{V,M}^{(\sigma,\sigma g)} Z_{M}(v,(\sigma g,\sigma)\tau).$$
Note that if $M\in  \mathscr{M}(\sigma g, \sigma)$ then $M$ is an irredcible $g$-twisted super $V$-module.
If $M$ is also irreducible $g$-twisted $V$-module then
$$Z_{M}(v,(\sigma g,\sigma),\tau)=\tr_M o(v) q^{L(0)-c/24}=\sum_{\lambda\in \Lambda_{M}}{\dim} W_{\l}Z_{M_{\l}}(v,\tau)$$
where $M_{\l}$ is an irreducible $V^G$-module by Theorem \ref{MTH1}.  If  $M$ is not an irreducible $g$-twisted $V$-module then
$M=N\oplus N\circ \sigma$ is a sum of two inequivalent irreducible $g$-twisted $V$-modules $N$ and $N\circ \sigma$
which are isomorphic $\bar g$-twisted $V_{\bar 0}$-modules and 
$$Z_{M}(v,(\sigma g,\sigma),\tau)=\tr_M o(v) q^{L(0)-c/24}=\sum_{\lambda\in \Lambda_{N}}2{\dim} W_{\l}Z_{N_{\l}}(v,\tau)$$
where $N_{\l}$ is an irreducible $V^G$-module by Theorem \ref{MTH1}.
Thus 
$$Z_{V^G}(v,-\frac{1}{\tau})=\frac{\tau^{\wt[v]}}{|G|}\sum_{j\in J, \lambda\in \Lambda_{M^j}}a_{j,\lambda}Z_{M^j_{\l}}(v,\tau).$$

Let $ {\cal M}_U$ be the set of inequivalent irreducible modules  for a regular vertx operator algebra $U$ of CFT type. Then $\{Z_W(v,\tau)|v\in U, W\in {\cal M}_U\}$ are linearly independent functions from $U\times \H\to \C$
and 
$$Z_{U}(v,-\frac{1}{\tau})=\tau^{\wt[v]}\sum_{W\in {\cal M}_{U}}S_{U,W}Z_W(v,\tau)$$
for $v\in U$ \cite{Z}. Moreover the coefficients $S_{V^G,W}\ne 0$  for all $W$ \cite{Hu}. Applying this result to $V^G$ implies that  any  irreducible $V^G$-module is isomorpic to one  in Theorem \ref{MTH1}. 
 \qed

 Th next result allows us to use various  results on the quantum dimensions  for vertex operator algebra
 $V^G$ \cite{DJX}.
 \begin{coro}  Let $V$ and $G$ be as in Theorem \ref{MTH1}. Then  the weight of every irreducible $V^G$-module is positive except $V^G$ itself.
 \end{coro}

\pf  It is clear from  \ref{MTH1} that $V^G$ is the only irreducible $V^G$-module whose weight  is $0$  and any other
irreducible $V^G$-module has a positive weight by Assumption (V2).
\qed
 
 From the discussion before, we know that  for any automorphism $g$ of $V$ of finite order, any irreducible $g\sigma^i$-twisted $V$-module $M$ restricts to an irreducible $\bar g$-twisted $V_{\bar 0}$-module if $M\circ \sigma\ne M,$
 and a direct sum of two inequivalent irreducible $\bar g$-twisted  $V_{\bar 0}$-modules if $M\circ \sigma=M$ where 
 $i=0,1.$ The following resuls tells us that every irreducible $\bar g$-twisted $V_{\bar 0}$-module is obtained in this way.

\begin{prop}\label{ptwisted}  Let $V$ be a regular, selfdual  \vosa of CFT type and $g$ an automorphismof $V$ of finite order.
Then every  irreducible $bar g$-twisted $V_{\bar 0}$-modules  occur in an irreduible  $g$ or $\sigma g$-twisted $V$-module.
\end{prop}
\pf 
Let $G$ be the group  generated by $g$ and $\sigma$ for $g\in G.$ Then $G$ is an abelian group and $V^G$ is a regular, selfdual vertx operator algebra of CFT type \cite{ABD, CM}. By Theorem 
 \ref{MTH1},  every irreudible $V^G$-module appears in an irreduible $h$-twisted $V$-module for some $h\in G.$ 
 Assume that there is an irreducible $\bar g$-twisted $V_{\bar 0}$-module $N$ which does not appear in any irreducibe 
 $g$-twisted  $V$-module.Then $N\circ \bar g\cong N$ and $N\circ\<\bar g\>=\{N\}.$ 
  Note that $V^G=V_{\bar 0}^{\<\bar g\>}.$ Let ${\cal A}$ be the set of irreducible $\bar h$-twisted $V_{\bar 0}$-modules  which appear  in irreducible $h$-twisted $V$-modules for all $h\in G.$ From Lemma \ref{mainlemma1} , we know  that ${\cal A}\circ \<\bar g\>={\cal A}.$ Since $N\notin {\cal A},$ any  irreducible $V^G$-module which occurs in $N$
  and any irreducible $V^G$-module which occurs in any $M\in {\cal A}$ are inequivalent. This is a contradiction as 
  every irreducible $V^G$-module is a $V^G$-submodule of some $M\in {\cal A}$ by Theorem  \ref{MTH1}. 
 \qed

We now discuss the modularity of $Z_M(v,(g,h),\tau)$ withou the assumption that $V^G$ is rational. 
\begin{prop}Let $V$ be a regular, selfdual  \vosa of CFT type and  $G$ a finite automorphism group of $V.$ Then there exists 
a congruence subgroup $A$ of $\Gamma$ such that for any $g,h\in G,$ $v\in V^G$ and $M\in \mathscr{M}(g, h)$
 $Z_M(v,(g,h),\tau)$
is a modular form of weight $\wt[v]$ over $A.$ Moreover, for each irreducible $V^G$-module $N$ appearing in an irreducible $g$-twisted $V$-module, $Z_N(v,\tau)$ is a modular form of weight $\wt[v]$ over $A.$
In particular, the character $\chi_{V^G}(\tau) $ of $V^G$ is  modular function over $A.$
\end{prop}
\pf In the case $V$ is a vertex operator algebra, the exactly same results were obtained in \cite[Lemma 4.3, Corollary 4.4]{DR}. The same proof 
works here.
\qed

\section{Quantum dimensions}
The quantum dimensions of the irreducible $g$-twisted $V$-modules and irreducible $V^G$-modules are computed in this section. We also present a super quantum Galois correspondence. The ideas and techniques used here come from \cite{DJX} and \cite{DRX1}.

Let $V$ be a vertex operator superalgebra as before and $G$ a finite automorphism group of $V.$ 
Let $g\in G$ and $M$ a $\sigma g$-twisted $V$-module. 
Recall $\chi_M(\tau)$ from Section 4,  We know that $\chi_V(\tau)$ and $\chi_M(\tau)$ are holomorphic functions on $\H.$ The quantum dimension of $M$ over $V$ is defined to be
$$\qdim_VM=\lim_{q\to 1^-}\frac{\ch_M(q)}{\ch_V(q)}=\lim_{y\to 0^+}\frac{\chi_M(iy)}{\chi_V(iy)}$$
where $q=e^{2\pi i\tau}$ and $y$ is real and positive.

We denote $M^s$ be the irreducible super $\sigma g$-twisted $V$-module
such that $M\subset M^s.$ Then $M^s=M$ if $M$ is $\sigma$-stable and $M^s=M\oplus M\circ \sigma$ if $M$ is not $\sigma$-stable. 
\begin{lem}\label{change} Let $M$ is a $\sigma g$-twisted $V$ -module which is not $\sigma$-stable.

(1) For $\l\in\Lambda_M$ set $M^s_{\l}=M_{\lambda}.$ Then $M^s=\oplus_{\lambda\in \Lambda_M}2W_{\lambda}	\otimes M^s_{\l}.$
	
 (2) For $\sigma h\in G_M.$ set
	$$Z_M(v,(g,h), \tau)=\tr_Mo(v)\phi(\sigma h)q^{L(0)-c/24}$$
for $v\in V_{\bar 0}.$ Then
	$$Z_{M^s}(v,(g,h), \tau)=2Z_M(v,(g,h), \tau)$$
	and $Z_{M^s}(v,(g,1), \tau)=0.$
	\end{lem}
\pf (1) is clear.
For (2) we only need to show that $Z_{M^e}(v,(g,1), \tau)=0$ for $v\in V_{\bar 0}.$ Note that
$M^s=N_1+N_2$ is a direct sum of two isomorphic irreducible $\overline{G}$-twisted $V$-modules such that 
$\sigma=(-1)^i$ on $N_i.$ It is immediate that $Z_{M^e}(v,(g,1), \tau)=0.$ 
\qed

The existence of the quantum dimension for a $\sigma g$-twisted $V$-module 
is given below in terms of the $S$-matrix and the proof is similar to that of \cite[Lemma 4.2]{DJX} 
by using the $S$-matrix given in Section 4. The proof requires Assumption (V2).
\begin{prop}\label{tqdim}
Let $M$ be an irreducible $\sigma g$-twisted $V$-module for $g\in G.$  Then $$\qdim_VM=\left\{\begin{array}{ll} \frac{S_{M,V}^{(g,\sigma)}}{S_{V,V}^{(\sigma,\sigma)}} & M=M^s\\
	\frac{1}{2}\qdim_VM^s=	\frac{1}{2}\frac{S_{M^s,V}^{(g,\sigma)}}{S_{V,V}^{(\sigma,\sigma)}} & M\ne M^s\end{array}\right.$$
\end{prop}
\pf Note that if $M$ is not $\sigma$-stable, then $\chi_M(\tau)=\frac{1}{2}\chi_{M^s}(\tau).$ So we can assume that
$M$ is $\sigma$-stable. 
A straightforward computation gives
\begin{equation*}\label{qdim S}
\begin{split}
\qdim_{V} M=&\lim_{y\to 0}\frac{\chi_{M}(iy)}{\chi_V(iy)}\\
=&\lim_{ y\to \infty}\frac{\chi_{M}(-\frac{1}{iy})}{\chi_{V}(-\frac{1}{iy})}\\
&\lim_{ y\to \infty}\frac{Z_M(\1, (g,\sigma), -\frac{1}{iy})}{Z_{V}(\1, (\sigma,\sigma),-\frac{1}{iy})}\\
=&\lim_{y \to \infty}\frac{\sum_{N\in  \mathscr{M}(\sigma , g^{-1}) }S_{M,N}^{(g,\sigma)} Z_{N}(\1,(\sigma,g^{-1}),iy)}
{\sum_{X\in \mathscr{M}(\sigma , \sigma)} S_{V,X}^{(\sigma,\sigma)} Z_X(\1,(\sigma,\sigma),iy)}\\
=&\lim_{y \to \infty}\frac{S_{M,V}^{(g,\sigma)}Z_{V}(\1,(\sigma,g^{-1}),iy)}{S_{V,V}^{(\sigma,\sigma)}Z_X(\1,(\sigma,\sigma),iy)}\\
=&\frac{S_{M,V}^{(g,\sigma)}}{S_{V,V}^{(\sigma,\sigma)}},
\end{split}
\end{equation*}
as desired.
\qed

Next we compute the quantum dimensions of irreducible $V^G$-modules. The ideas and techniques come from
the proof of \cite[Theorem 4.4]{DRX1}.
Recall from \cite[Lemma 4,3]{DRX1} that for any  $\chi_1, \chi_2$ the two irreducible characters of  $\C^{\a_M}[G_M]$,  
$\frac{1}{|G_M|}\sum_{a\in G_M}\chi_1(\hat a)\overline{\chi_2(\hat a)}=\delta_{\chi_1,\chi_2}$ where $\overline{\chi_2(\hat a)}$ is the complex conjugation of $ \chi_2(\hat a).$ We remark that 
$\bar a$ was used in \cite{DRX1} instead of $\hat a.$  We use $\bar a$ for the element in $\overline{G}$  in this paper.

\begin{thm}\label{MTH} Let $M$ be an irreducible $\sigma g$-twisted $V$-module and $\lambda\in \Lambda_{G_M}.$  Then 
$$S_{M_\l,V^G}=\left\{\begin{array}{ll}\frac{\dim W_{\l}}{|G_M|} S_{M,V}^{(g,\sigma)}& M=M^s\\
\frac{\dim W_{\l}}{2|G_M|}S_{M^s,V}^{(g,\sigma)}& M\ne M^s\end{array}\right.$$
and $\qdim_{V^G}M_{\l}=[G:G_M]\dim{W_\l}\qdim_VM.$ In particular, $\qdim_{V^G}V_{\l}=\dim{W_\l}.$
Moreover,  $\qdim_{V^G}M=|G|\qdim_VM .$
\end{thm}
\pf 
(1) Assume that $M$ is $\sigma$-stable. Recall the decomposition
$$M=\oplus_{\mu\in\Lambda_{M}}W_{\mu}\otimes M_{\mu}.$$
For $v\in V^G$ using orthogonal relation of two irreducible characters of $\C^{\a_M}[G_M]$ we see that
 \begin{equation*}
\begin{split}
\frac{1}{|G_M|}\sum_{h \in G_M}Z_M(v, (g,h),\tau)\overline{\lambda(\widehat{\sigma h})}
=&\frac{1}{|G_M|}\sum_{h \in G_M}\sum_{\mu\in {\Lambda_M}}\mu(\widehat{\sigma h })\overline{\lambda(\widehat{\sigma h})}Z_{M_\mu}(v, \tau)\\
=&Z_{M_\l}(v, \tau)
 \end{split}
\end{equation*} 
where $\overline{\lambda(\widehat{\sigma h})}$ is the complex conjugation of $\lambda(\widehat{\sigma h}).$
By Theorem \ref{minvariance} we have
\begin{equation*}
\begin{split}
Z_{M_\l}(v, -1/\tau)&=\frac{1}{|G_M|}\tau^{\wt[v]}\sum_{h\in G_M}\sum_{N\in  \mathscr{M} (h,g^{-1})}S_{M,N}^{(g,h)}Z_N(v,(h,g^{-1}),\tau)\overline{\lambda(\widehat{\sigma
	h})}.
\end{split}
\end{equation*}
Note that if  $h=\sigma$ and $\overline{\lambda(\widehat{\sigma
		h})}=\dim W_{\l}.$  So
the coefficient of $Z_{V^G}(v,\tau)$ in $Z_{M_\l}(v, -1/\tau)$ is $\tau^{\wt[v]}\frac{\dim W_\l}{|G_M|}S_{M,V}^{(g,\sigma)}.$

On the other hand,
$$Z_{M_\l}(v, -1/\tau)=\tau^{\wt[v]}\sum_{W\in {\cal M}_{V^G}}S_{M_\l,W}Z_W(v,\tau)$$
where ${\cal M}_{V^G}$ is the set of inequivalent irreducible $V^G$-modules and 
$(S_{X,Y})_{X,Y\in  {\cal M}_{V^G}}$ is the $S$-matrix of $V^G.$
We see that $S_{M_\l,V^G}=\frac{\dim W_{\l}}{|G_M|}S_{M,V}^{(g,\sigma)}.$ 
 In particular,  $S_{V^G,V^G}= \frac{1}{|G|}S_{V,V}^{(\sigma,\sigma)}.$
It follows from Proposition \ref{tqdim} that
$$\qdim_{V^G}M_{\l}=\frac{S_{M_\l,V^G}}{S_{V^G,V^G}}
=\frac{\frac{\dim W_{\l}}{|G_M|}S_{M,V}^{(g,\sigma)}}{\frac{S_{V,V}^{(\sigma,\sigma)}}{|G|}}
=[G:G_M]\dim{W_\l}\qdim_VM.$$
The equality $\qdim_{V^G}M=|G|\qdim_VM $ follows from the fact that $|G_M|=\sum_{\lambda\in\Lambda_{M}}(\dim W_\l)^2.$

(2) Assume that $M$ is not $\sigma$-stable.   Using Lemma \ref{change} yields
\begin{equation*}
\begin{split}
Z_{M_\l}(v, \tau)=&\frac{1}{|G_M|}\sum_{h \in G_M}\sum_{\mu\in {\Lambda_M}}\mu(\widehat{\sigma h })\overline{\lambda(\widehat{\sigma h})}Z_{M_\mu}(v, \tau)\\
=&\frac{1}{2}
\frac{1}{|G_M|}\sum_{h \in G_M}Z_{M^s}(v, (g,h),\tau)\overline{\lambda(\widehat{\sigma h})}.
\end{split}
\end{equation*} 
From the proof of (1) we see that 
\begin{equation*}
\begin{split}
Z_{M_\l}(v, -1/\tau)&=\frac{1}{2|G_M|}\tau^{\wt[v]}\sum_{h\in G_M}\sum_{N\in  \mathscr{M} (h,g^{-1})}S_{M^s,N}^{(g,h)}Z_N(v,(h,g^{-1}),\tau)\overline{\lambda(\widehat{\sigma
		h})}.
\end{split}
\end{equation*}
This implies that $S_{M_\l,V^G}=\frac{\dim W_{\l}}{2|G_M|}S_{M^s,V}^{(g,\sigma)}$ and $\qdim_{V^G}M_{\l}=[G:G_M]\dim W_{\l}\qdim_{V}M$ by using Proposition \ref{tqdim}.
The equality $\qdim_{V^G}M=|G|\qdim_VM $ again follows  from the fact that $|G_M|=\sum_{\lambda\in\Lambda_{M}}(\dim W_\l)^2.$
\qed

\begin{lem}\label{positive} Let $U$ be a selfdual, regular vertex operator algebra of CFT type and $(S_{X,Y})_{X,Y\in  {\cal M}_U}$ be the $S$-matrix of $U.$ Then $S_{X,U}>0$ for any $X\in {\cal M}_U.$
\end{lem}

\pf From \cite{DJX} we know that $\frac{S_{X,U}}{S_{U,U}}>0$ for all $X\in {\cal M}_U$ and $S_{U,U}$ is a real number. This implies that all $S_{X,U}$ are positive or negative.  Since $\frac{-1}{i}=i$ we see that
$$\chi_V(i)=\sum_{X\in {\cal M}_U}S_{U,X}\chi_X(i),$$
Note that $\chi_X(i)>0$ for all $X.$ This forces $S_{U,X}=S_{X,U}$ to be positive for all $X.$
\qed

The following corollary is immediate from Theorem \ref{MTH} and Lemma \ref{positive}.
\begin{coro} Let $V$ and $G$ be as before. Then for any irreducible $g$-twisted super $V$-module $M,$ 
$S_{M,V}^{(g,\sigma)}$ is positive.
\end{coro}

The next result tells us the tensor product of two irreducible $V^G$-modules in $V$ is determined by the 
the tensor product of corresponding irreducible $G$-modules.
\begin{thm} For $\l,\mu\in \irr(G),$ 
	$$V_{\l}\boxtimes V_{\mu}\cong \sum_{\gamma\in \irr(G)}{\rm Hom}_{G}(W_\l\otimes W_\mu, W_{\gamma}) V_{\gamma}$$
\end{thm}

\pf The proof of \cite[Corollary 5.7]{DRX1} works here (aslo see \cite{T}).
\qed
	
We are now ready to discuss the super quantum Galois theory. For any VOA extension $V\supset U$ we can define the Galois group $\Gal(V/U)=\{g\in \Aut(V)\left|\ g|_{U}=\id\right \}$ as in the classical field theory. 
We also define the
index $[V:U]=\qdim_UV.$ $V$ is called a Galois extension of $U$ if $|\Gal(V/U)|=[V:U].$
\begin{thm}\label{Galois} The map $H\to V^H$ gives a bijection between the set of subgroups of $G$ and the set
	of vertex operator super subalgebras of $V$  satisfying the following:
	
	(1) $\sigma \in H$  if and only if $V^H\subset V_{\bar 0},$
	
	(2) $[V:V^H]=o(H)$ and   $[V^H:V^G]=[G:H]$,
	
	(3) $H\lhd G$ if and only if $V^H$ is a Galois extension of $V^G$. In this case $\Gal(V^H/V^G)$ is isomorphic $G/H.$
\end{thm} 

\pf (1) It is clear that $V^H\subset V_{\bar 0}$ if $\sigma \in H.$ Now assume that $V^H\subset V_{\bar 0}$ and 
$\sigma \notin H.$ Then $V_{\bar i} =\sum_{\l\in \irr(H)}W_{\l}\otimes (V_{\bar i})_{\l}$ for $i=0,1$  and  $(V_{\bar i})_{1}\ne 0$ where $1$ is the trivial character of $H$ by \cite[Theorem 5.4]{DY}.
So $V^H=V_{\bar 0}^H+V_{\bar 1}^H=(V_{\bar 0})_1+(V_{\bar 1})_1$ is not a subspace of $V_{\bar 0}$ and this is a contradiction.

(2) By Theorem \ref{MTH}, $[V:V^H]=\qdim_{V^H}V=o(H)\qdim_VV=o(H)$ and $[V^H:V^G]=\qdim_{V^G}V^H=
\frac{\qdim_{V^G}V}{\qdim_{V^H}V}=[G:H].$ 

(3) The proof is similar to that of \cite[Theorem 6.10 (2)]{DJX}.
	
Finally we show  the map $H\mapsto V^H$ is a bijection.  
For convienience, we let ${\cal G}_1=\{H<G|\sigma \in H\}, {\cal G}_2=\{H<G|\sigma \notin H\}$
and ${\cal V}_1$ be the set of vertex operator subalgebras of $V_{\bar 0}$ which contain $V^G,$
${\cal V}_2$  be set of vertex operator super subalgebras of $V$ which are not in ${\cal V}_1.$ It is good enough to show 
that the map $H\mapsto V^H$ gives a bijection from ${\cal G}_s$ to ${\cal V}_s$ for $s=1,2.$
By \cite[Theorem 6.10]{DJX}, $A\mapsto V_{\bar 0}^{A}$ 
is a bijection from the set of subgroups of $\overline{G}$ and vertex operator subalgebras of $V_{\bar 0}$ containing $V_{\bar 0}^{\overline{G}}=V^G.$ 

For $s=1,$ we notice that $\overline{G}$ is isomorphic to $G/\<\sigma\>.$ For any subgorup $A<\overline{G}$ there exists
$H\in {\cal G}_1$ such that  $A=\overline{H}$ and $V^A=V^H.$ If $H_1,H_2\in {\cal G}_1$ such that $\overline {H_1}=\overline {H_2},$
Then for any $h_1\in H_1,$ $\bar h_1=\bar h_2$ for some $h_2\in H_2.$ This implies that $h_1=h_2\sigma^s\in H_2$ for some $s$ and  $H_1=H_2.$ So $H\mapsto V^H$ gives a bijection from ${\cal G}_1$ to ${\cal V}_1.$

For $s=2$ we know from (1) that if $H\in {\cal G}_2$ then $V^H\in {\cal V}_2.$ Assume that $H_1,H_2\in {\cal G}_2$ such that $V^{H_1}=V^{H_2}.$ Note that $V_{\bar 0}^{\overline{H}_1}=V^{H_1}_{\bar 0}=V_{\bar 0}^{H_2}=V_{\bar 0}^{\overline{H_2}}$ where $\overline{H_i}$ is the image of $H_i$ in $\overline{G}.$ From the discussion for $s=1,$ we have  $\overline {H_1}=\overline {H_2}.$  If $H_1\ne H_2$
then there exists $h_i\in H_i$ such that $h_1=h_2\sigma.$  From $V^{H_1}\cap V_{\bar 1}=V^{H^2}\cap  V_{\bar 1}\ne 0$  by (1), we see that $h_1=1$ and $-1$ on $V^{H_1}\cap V_{\bar 1},$ a contradiction. So $H\mapsto V^H$
gives an injection from ${\cal G}_2$ to ${\cal V}_2.$ It remains to show the map is onto. Let $U$ be a vertex operator super subalgebra of $V$ such that $U$ contains $V^G$ and $U_{\bar 1}\ne 0.$ The problem is that we even do not know if $U$ is a simple. Then there exists $H<{\cal G}_1$ such  that $U_{\bar 0}=V^H$ is a simple vertex operator algebra.  We claim that there is a subgroup $A$ of $H$ of index 2 such that $U=V^A.$ A proof of this claim is divided into three steps:

1) From the decomposition $V=\oplus_{\lambda\in\irr(H)}W_{\lambda}\otimes V_{\lambda}$ we know that $\{V_{\l}|\lambda\in \irr(H)\}$  are inequivalent irreducible $V^H$-modules. Recall from Section 5 that $V^{\l}=W_{\l}\otimes V_{\l}$ is the sum of simple $H$-modules of $V$ isomorphic to $W_{\l}.$ Since $hY(u,z)v=Y(hu,z)hv$ for any $u\in V^{\l},$ $v\in V^{\mu}$ we see that 
$$V^{\l}\cdot V^{\mu}=\<u_nv|u\in V^{\l}, v\in V^{\mu}, n\in\Z\><\sum_{\gamma\in \l\otimes \mu}V^{\gamma}.$$
Set $U^{\l}=U\cap V^{\l}.$ Then $U_{\bar 1}=\sum_{\l\in \irr(H)_1}U^{\l}$ where $\irr(H)_1$ consists of
irreducible characters $\l$ such that $\sigma=-1$ on $W_{\l}.$ By \cite[Lemma 3.2]{DM1}, we know that
$$0\ne U^\l\cdot U^{\mu}<\sum_{\gamma\in \l\otimes \mu}U^{\gamma}$$
 if $U^\l, U^\mu\ne 0.$  Assume that $U^{\l}, U^{\mu}$
are nonzero for two distict $\l,\mu\in \irr(H)_1.$ Then either $\l$ is not equal to its dual $\l^*$ or $\l\ne \mu^*.$ 
If $\l\ne \l^*$ then $$0\ne U^\l\cdot U^{\l}<U_{\bar 1}\cdot U_{\bar 1}=V^H.$$
 But $\l\otimes \l$ does not contain
the trivial character $1$ and this is a contradiction. Similarly we have a contradiction if $\l\ne \mu^*.$ This imples
that $U_{\bar 1}=U^{\l}$ for some $\l\in \irr(H)_1$ satisfying $\l=\l^*.$

2) We prove in this step that $U_{\bar 1}=V_{\l}=V^\l$ and $\l$ is a linear character.  Write $U^{\l}=X_{\l}\otimes V_{\l}$ where
$X_\l$ is a subspace of $W_{\l}.$ Recall from the proof of \cite[Proposition 2.1]{DM2} that
 regarding $W_{\l}$ as a homogeneous subspace of $V^{\l}$ as $H$-modules, then there exists integer $n$
such that 
$Y_n=\<\sum_{m\geq n}u_mv|u,v\in W_{\l}\>$ is isomorphc to $W_{\l}\otimes W_{\l}$ as $H$-modules by identifying 
	$\sum_{m\geq n}u_mv$ with $u\otimes v.$ Recall from group theory that the trivial module $\C$ has multiplicity 
	1 in $W_{\l}\otimes W_{\l}$ spanned by $\sum_{i=1}^{d}u^i\otimes u^{i*}$ where $\{u^1,...,u^d\}$ is a basis of $W_\l$
	and $\{u^{1*},...,u^{d*}\}$ is a dual basis of $W_{\l^*}=W_{\l}.$ Note that $U^{\l}$ is a direct sum of infinitely many
	copies of $X_{\l}.$ For convenience, we can also assume that $X_{\l} $ is a subspace of $W_{\l}.$ Then 
	$Z_n=\<\sum_{m\geq n}u_mv|u,v\in X_{\l}\>$ is isomorhic to $X_l\otimes X_{\l}$ from the discussion above.
	On the other hand $Z_n<U_{\bar 0}=V^H.$ This forces $X_{\l}=W_{\l}$ and $\l$ is a linear character.
	
3) Let $A=\ker \l.$ Since $\l^2=1,$ $\l(h)=\pm 1$ for all $h\in H.$ Also $\l(\sigma)=-1.$ So $A$ is a subgroup of $H$ 
of index 2. It is clear now that $U=V^A.$ That is , the map $H\mapsto V^H$ is onto from
${\cal G}_2\to {\cal V}_2.$ the proof is complete.
\qed

\section{Global dimensions}

We study the global dimensions of $V$ and $V^G$ in this section.  We define the global dimension of $V$ $$\glob(V)=\sum_{M\in {\cal S}(1)}(\qdim_VM)^2$$
where ${\cal S}(g)$ is the set  of inequivalent irreducible $g$-twisted $V$-module. 
In the case $V=V_{\bar 0}$ is a vertex operator algebra, this is exactly the global dimension of $V$ defined in \cite{DJX} and It is proved there  that $\glob(V)=\frac{1}{S_{V,V}^2}.$
\begin{lem}\label{lglo1}  We have the following relation: 
	
	(1) $\glob(V)=\left(\frac{1}{S_{V, V}^{(\sigma, \sigma)}}\right)^{2}.$
	
	(2) $\glob(V^G)=|G|^2\glob(V)=|G|^2\sum_{j\in J}\frac{1}{|G_{M^j}|}(\qdim_VM^j)^2.$
\end{lem}
\pf 
(1)	 was given in the proof of Proposition 7.2 in \cite{DNR1}.

(2) Using Theorem \ref{MTH} and (1) yields 
$$\glob(V^G)=\frac{1}{(S_{V^G, V^G})^2}=\frac{|G|^2}{\left(S_{V, V}^{(\sigma, \sigma)}\right)^{2}}=|G|^2\glob(V).$$
A direct calculation  using Theorem  \ref{MTH}  gives 
$$\glob(V^G)=\sum_{j\in J, \lambda\in \Lambda_{M^j}}(\qdim_{V^G}M^j_{\l})^2=\sum_{j\in J} \frac{|G|^2}{|G_{M^j}|}(\qdim_{V}M^j)^2.$$
In particular, we have a new formula 
$$\glob(V)=\sum_{j\in J} \frac{1}{|G_{M^j}|}(\qdim_{V}M^j)^2=\frac{1}{|G|}\sum_{M\in {\cal S}}(\qdim_{V}M)^2$$
where we have used the fact that $\sum_{M\in M^j\circ G}(\qdim_VM)^2=[G:G_{M^j}](\qdim_VM^j)^2.$
\qed

We reamrk that Lemma \ref{lglo1} (2)  for vertex operator agebra was established previously in \cite{DRX1}.

We investigate more on the global dimension $\glob(V)$ in terms of twisted modules. In particular, we show that 
for any automorphism $g$ of $V$ of finite order,
$$\glob(V)=\sum_{M\in {\cal S}(g)}(\qdim_VM)^2.$$
The main idea is to prove that the sum
of the squares of the quantum dimensions of irreducible $V^G$-modules appearing in irreducible $h$-twisted modules is $\frac{1}{o(G)}\glob(V^G)$ for any $h\in G$ where $G$ is the subgroup of $\Aut(V)$ generated by $g$ and $\sigma.$


 
 Note that  $V$ has a decomposition $V=\oplus_{\chi \in \irr(G)}W_{\chi}\otimes V_{\chi}.$
 Since $G$ is an abelian group, each $W_{\chi}$ is one-dimensional. So $V_{\chi}$ is a simple current \cite{DJX}.
 \begin{lem}\label{eigen} Let $h\in G$   and $M$ an irreducible $\sigma h$-twisted $V$-module, $\chi\in \irr(G), \l\in G_M.$ Then
 	$\frac{S_{V_{\chi}, M_{\lambda}}}{S_{V^G, M_{\lambda}}}=\overline{\chi(\sigma h)}.$ 
 \end{lem}
 \pf  Note that for $v\in V^G$ 
 $$Z_{V_{\chi}}(v, \tau)=\frac{1}{|G|}\sum_{x\in G}Z_V(v, (\sigma,x),\tau)\overline{\chi(\sigma x)},$$
 $$Z_{V^G}(v, \tau)=\frac{1}{|G|}\sum_{x\in G}Z_V(v, (\sigma,x),\tau).$$
 Thus,
 $$
 Z_{V_{\chi}}(v, -1/\tau)=\frac{\tau^{\wt[v]}}{|G|}\sum_{x\in G}\sum_{N\in \mathscr{M}(x, \sigma)}S_{V,N}^{(\sigma, x)}Z_N(v,(x,\sigma)\tau)\overline{\chi(\sigma h)},$$
 $$Z_{V^G}(v, -1/\tau)=\frac{\tau^{\wt[v]}}{|G|}\sum_{x\in G}\sum_{N\in \mathscr{M}(x, \sigma)}S_{V,N}^{(\sigma, x)}Z_N(v,(x,\sigma)\tau).$$
 Note that $\mathscr{M}( x, \sigma)\circ G=\mathscr{M}( x, \sigma)$ and $N$ is a completely reducible $V^G$-module 
 for $N\in \mathscr{M}( x, \sigma).$ The result follows.
 \qed

 Motivated by a similar resut  for a vertex operator algebra from \cite{DRX1} we have the following result with a similar proof.
 \begin{thm}\label{class}  Let $h\in G.$ Then
 	$$\sum_{X\in {\cal M}_{V^G}^h }(\qdim_{V^G}X)^2=\frac{1}{|G|} \glob(V^G)$$
 	where  ${\cal M}_{V^G}^h$ is the set of inequivalent irreducible $V^G$-modules appearing in the irreducible $h$-twisted $V$-modules.
 \end{thm}
\pf   Since
$$\sum_{X\in {\cal M}_{V^G}^h}(\qdim_{V^G}X)^2=\frac{1}{S_{V^G,V^G}^2}\sum_XS_{V^G,X}^2=\glob(V^G)\sum_XS_{V^G,X}^2$$
it suffices to show that
$$\sum_{X}S_{V^G,X}^2=\frac{1}{|G|}.$$
 	Set $x_h=\sum_XS_{V^G,X}^2.$
 	Using  Lemma \ref{positive}, Theorem \ref{MTH} , and the unitarity of $S$-matix  gives the  orthogonal relation
 	$$\sum_{Z \in {\cal M}_{V^G}}S_{V_{\chi},Z}S_{V^G,Z}=\delta_{\chi,1}$$
 	for any $\chi$ where $1$ is the trivial character of $G.$ 
 	
 	From Theorem \ref{MTH1}, ${\cal M}_{V^G}$ is a disjoint  union of ${\cal M}_{V^G}^{h}$
 	for $h\in G.$  Using Lemma \ref{eigen} gives a linear system
 	$$\sum_{h\in G}x_h\overline{\chi(h)}=\sum_{h\in G}\sum_{X\in {\cal M}_{V^G}^h}S_{V^G,X}^2
 	\frac{S_{V_{\chi},X}}{S_{V^G,X}}=\delta_{\chi ,1}$$
 	for $\chi\in\irr(G)$ wth coefficient matrix $A=(\chi(h))_{h\in G, \chi\in\irr(G)}.$  Since $\sum_{h\in G}\chi_1(h)\overline{\chi_2(h)}=\delta_{\chi_1,\chi_2}|G|$ we see that $A$ is nondegenerate. 
 	So the linear system has a unique solution $(x_h)_{h\in G}.$ It is easy to see that
 	$x_h=\frac{1}{|G|}$ is a solution. The proof is complete.
 \qed

The following result is an extension of the same result for vertex operator algebra \cite[Theorem 5.3]{DRX1} to \vosa.
\begin{thm} \label{main} Let $V$ be a selfdual, regular \vosa of CFT type, $g$ an automorphism of $V$ of finite order 
 and the  the weight of any irreducible $h$-twisted $V$-module for $h\in G$ is positive except $V$ itself  where $G$ is  the group generated by $g$ and $\sigma.$ Then
 $$ \glob(V)=\sum_{M\in {\cal S}(g)}(\qdim_VM)^2.$$
\end{thm}
\pf   Since $G$ is abelian group, ${\cal S}(g)\circ G={\cal S}(g).$ Let $J_g=\{j\in J| M^j\in {\cal S}(g)\}.$
Applying Theorems \ref{class}, \ref{MTH} and Lemma \ref{lglo1}  we see that
\begin{equation*}
\begin{split}
\glob(V)=&\frac{1}{|G|^2}\glob(V^G)\\
=&\frac{1}{|G|}\sum_{j\in J_g, \lambda\in \Lambda_{M^j}}(\qdim_{V^G}M^j_\l)^2\\
=&\frac{1}{|G|}\sum_{j\in J_g, \lambda\in \Lambda_{M^j}}[G:G_{M^j}]^2(\dim W_{\l})^2(\qdim_{V}M^j)^2\\
=&\sum_{j\in J_g}[G:G_{M^j}](\qdim_VM^j)^2\\
=&\sum_{M\in {\cal S}(g)}(\qdim_{V}M)^2,
\end{split}
\end{equation*}
as desired.
\qed

We remark that in the case $g=\sigma,$ these results have been established in \cite{DNR1} previously.  

The following result generalizes Theorem \ref{class} from cyclic group to any finite group.
\begin{coro} Let $G$ be a finite group. Then for any $g\in G,$
$$\sum_{X\in {\cal M}_{V^G}^g}(\qdim_{V^G}X)^2=\frac{1}{|C_G(g)|} \glob(V^G)$$
where $ {\cal M}_{V^G}^g$ is the set of  inequivalent irreducible $V^G$-modules appearing in the irreducible $g$-twisted $V$-modules.
\end{coro}
\pf  Based on Theorems \ref{class} and \ref{main}, the same proof of \cite[Corollary 5.4]{DRX1} works here. 
\qed

Recall that a \vosa $V$ is called holomorphic if it is rational and $V$ is the only irreducible module (up to isomorphism)
for itself.  If $V$ is a vertex operaqtor algebra there is a unique irreducible $g$-twisted $V$-module $V(g)$ such that $\qdim_V(g)=1$ for any automorphism 
$g$ of $V$ of finite order \cite{DLM7, DRX1}. But if $V$ is not a vertex operator algebra $(V_{\bar 1}\ne 0),$ this is not always ture.  From Theorem \ref{minvariance}, $V$ has a unique irreducible $g$-twisted super $V$-module $V(g).$
Then $\qdim_VV(g)=1$ if  $V(g)$ is  an unique irreducible $g$-twisted $V$-module, and  $\qdim_VV(g)=\sqrt{2}$ if $V(g)$ is a direct sum of two inequivalent irreducible $g$-twisted $V$-modules by Theorem \ref{MTH1}.

\section{$S$-matrix}

By Theorem \ref{MTH1}, $\left\{M_{\lambda}^{j} \mid j \in J, \lambda \in \Lambda_{M^j}\right\}$ gives a complete list of inequivalent 
irreducible $V^G$-modules. We determine the $S$-matrix  $\left(S_{M_{\lambda}^{i}, M_{\mu}^{j}}\right)_{(i, \lambda),(j, \mu)}$ of $V^G$ in this section.    The proof of Theorem \ref{MTH} essentially gives the $S$-matrix, but we want to find an explict formula of $S_{M^i_{\l}, M^j_{\mu}}$  in terms of the $S$-matrix given in Theorem \ref{minvariance} following the ideas and techniques in \cite{DRX2}.

Note that $Z_{M^i_{\l}}(v,\tau)=\tr_{M^i_{\l}}o(v)q^{L(0)-c/24}$ for $v\in V^G.$ 
Using Theorem \ref{minvariance}, we have
$$Z_{M^i_{\l}}(v,\tau+1)=e^{2\pi i(-c/24+h_{i,\l})}Z_{M^i_{\l}}(v,\tau),$$
$$
Z_{M_{\lambda}^{i}}(v,-\frac{1}{\tau})=\tau^{\mathrm{wt}[v]} \sum_{j \in J, \mu \in \Lambda_{M^j}} S_{M_{\lambda}^{i}, M_{\mu}^{j}} Z_{M_{\mu}^{j}}(v, \tau)
$$
where $c$ is the central charge of $V,$ $h_{i,\l}$ is the conformal weight of $M^io_{\l}$.
The $T$ matrix of $V^G$ is  diagonal unitary matrix with $T_{M^i_{\l},M^i_{\l}}=e^{2\pi i(-c/24+h_{i,\l})}$.
$\left(S_{M_{\lambda}^{i}, M_{\mu}^{j}}\right)_{(i, \lambda),(j, \mu)}$ is the $S$ matrix of $V^G$,
which is the key to understand the action of $\Gamma$.

Let $M$ be an irreducible $\sigma g$-twisted $V$-module. Following \cite[Lemma 4.3]{DNR1}  we take $\alpha_{M\circ k}$ such that $\alpha_{M\circ k}(k^{-1}ak,k^{-1}bk)=\alpha_{M}(a,b)$
for $a,b\in G_M$ and $k\in G.$ 
Then $\C^{\a_M}[G_M]$ and  $\C^{\a_{M\circ k}}[G_{M\circ k}]$ are isomorphic by sending $\hat a\in \C^{\a_M}[G_M]$ to $\widehat{k^{-1}ak}.$ For any $\lambda\in \Lambda_{G_M}$ we define  $\lambda\circ k\in  \Lambda_{G_{M\circ k}}$ such that $ (\lambda\circ k)(\widehat{k^{-1}ak})=
\lambda(\hat{a})$  for  $a\in G_M.$
It is easy to check that $Z_{(M\circ k)_{\lambda\circ k}}(v,\tau)=Z_{M_{\lambda}}(v,\tau).$ Moreover, if $M\in \mathscr{M}(g,h)$ then  $Z_{M\circ k}(v, (k^{-1}gk,k^{-1}hk),\tau)=Z_M(v, (g,h),\tau).$

\begin{lem}\label{SMN}
	Let $M\in \mathscr{M}(g,h)$ and $N\in \mathscr{M}(h, g^{-1})$. 
	Then  for any $k\in  G$, $S_{M,N}^{(g,h)}=S_{M\circ k,N\circ k}^{(k^{-1}gk,k^{-1}hk)}$.
\end{lem}
\pf
	Using  Theorem \ref{minvariance} gives
	$$Z_M(v, (g,h), -1/\tau)=\tau^{\wt[v]}\sum_{N\in \mathscr{M}(h, g^{-1})}S_{M,N}^{(g,h)} Z_{N}(v,(h,g^{-1}),\tau).$$
 On the other hand,
	\begin{equation*}
	\begin{split}
	Z_M(v, (g,h), -1/\tau)=&Z_{M\circ k}(v, (k^{-1}gk,k^{-1}hk),\tau)\\
	=&\tau^{\wt[v]}\sum_{N\in \mathscr{M}(h, g^{-1}) }S_{M\circ k,N\circ k}^{(k^{-1}gk,k^{-1}hk)}Z_{N}(v,(h,g^{-1}),\tau).
	\end{split}
	\end{equation*}
	Since 
	$\{Z_{N}(v,(h,g^{-1}),\tau)  \mid N\in \mathscr{M}(h, g^{-1}) \}$ are linearly independent functions on
	$V\times \H$  \cite{DZ1}, the result follows.
\qed

For $i,j\in J$, let $M^j$ be a $\sigma g_j$-twisted $V$-module and $C_{i,j}$ an least subset of $G$ such that $\{M^j\circ k \mid k\in C_{i,j}\}=M^j\circ G\cap (\cup_{h\in G_{M^i}} {\cal S}(\sigma h, \sigma g_i^{-1}))$ \cite{DRX2} 
where  ${\cal S}(\sigma h, \sigma g_i^{-1})=\{N\in {\cal S}(\sigma h)|N\circ \sigma g_i^{-1}=N\}.$
Note that $\{(M^j)^s\circ k \mid k\in C_{i,j}\}=(M^j)^s\circ G\cap (\cup_{h\in G_{M^i}} \mathscr{M}(\sigma h, \sigma g_i^{-1})).$ It is possible that $C_{i,j}$ is an empty set.
Assume that $0\in J$ and $M^0=V.$ 
In this case, $G_{M^i}=G$ and $\{M^j\circ k \mid k\in C_{0,j}\}=M^j\circ G$ and  $|C_{0,j}|=[G,G_{M^j}]$.

Recall Lemma \ref{change}.  We will let $\Lambda_{M^s}=\Lambda_M$ when $M\ne M^s$ in the following discussion.
\begin{thm}\label{ZVG}
	Let $i,j\in J$ and $\lambda\in\Lambda_{G_{M^i}}$, $\mu\in \Lambda_{G_{M^j}}.$ Set $a_i=1$ if $M^i=(M^i)^s$
	and $a_i=2$ if $M^i\ne (M^i)^s.$
	Then
	\begin{equation*}
	S_{M^i_\l,M^j_\mu}=\left\{\begin{array}{ll}
	\frac{a_ja_i^{-1}}{|G_{M^i}|}\sum_{k\in C_{i,j}}S_{(M^i)^s,(M^j)^s\circ k}^{(g_i,k^{-1}g_j k)}\overline{\lambda(\widehat{k^{-1}g_jk})}\mu(\widehat{kg_i^{-1}k^{-1}}) & C_{i,j} \ne \emptyset\\
	\ &\\
	0 &   C_{i,j}=\emptyset
	\end{array}\right.
	\end{equation*}
	In particular,  for $i=0$ and $j\in J$, 
	$S_{M^0_\l,M^j_\mu}=\frac{a_j}{|G_{M^j}|}S_{M^0,M^j}^{(1,g_j)}\lambda(g_j^{-1})\dim W_{\mu}.$
	\end{thm}
\pf From the proof of Theorem \ref{MTH} we have 
	$$
	Z_{M^i_\l}(v, \tau)=\frac{1}{a_i|G_{M^i}|}\sum_{h\in G_{M^i}}Z_{(M^i)^s}(v, (g_i,h),\tau)\overline{\lambda(\widehat{\sigma h})}.
	$$
	and 
	\begin{equation*}
	\begin{split}
Z_{M^i_\l}(v, -1/\tau)=&\frac{\tau^{\wt[v]}}{a_i|G_{M^i}|}\sum_{h\in G_{M^i}}\sum_{N\in \mathscr{M} (h,g^{-1})}S_{(M^i)^s,N}^{(g_i,h)}Z_N(v,(h,g_i^{-1}),\tau)\overline{\lambda(\widehat{\sigma h})}\\
=&\frac{\tau^{\wt[v]}}{a_i|G_{M^i}|}\sum_{j\in J}\sum_{k\in C_{i,j}}S_{(M^i)^s,(M^j)^s\circ k}^{(g_i,k^{-1}g_j k)}
Z_{(M^j)^s\circ k}(v,(k^{-1}g_jk, ,g_i^{-1}),\tau)\overline{\lambda(\widehat{\sigma h})}\\
=&\frac{\tau^{\wt[v]}a_ja_i^{-1}}{|G_{M^i}|}\sum_{j\in J}\sum_{k\in C_{i,j}}S_{(M^i)^s,(M^j)^s\circ k}^{(g_i,k^{-1}g_j k)}\
Z_{M^j\circ k}(v,(k^{-1}g_jk, ,g_i^{-1}),\tau)\overline{\lambda(\widehat{\sigma h})}\\
\end{split}
\end{equation*}
where we have identified  $N, h$ with  $(M^j)^s\circ k, k^{-1}g_jk$ respectively.
So
	\begin{equation*}
	\begin{split}
	Z_{M^i_\l}(v, -1/\tau)
	&=\frac{\tau^{\wt[v]}a_ja_i^{-1}}{|G_{M^i}|}\sum_{j\in J}\sum_{k\in C_{i,j}}\sum_{\mu\in \Lambda_{M^j\circ k}}S_{(M^i)^s,(M^j)^s\circ k}^{(g_i,k^{-1}g_j k)}\overline{{\lambda(\widehat{k^{-1}g_jk})}}{\mu(\widehat {g_i^{-1}})}Z_{(M^j\circ k)_{\mu}}(v,\tau)\\
	&=\frac{\tau^{\wt[v]}a_ja_i^{-1}}{|G_{M^i}|}\sum_{j\in J}\sum_{k\in C_{i,j}}\sum_{\mu\in \Lambda_{{M^j}}}S_{M^i,M^j\circ k}^{(g_i,k^{-1}g_j k)}\overline{{\lambda(\widehat{k^{-1}g_jk})}}\mu(\widehat {kg_i^{-1}k^{-1}})Z_{M^j_{\mu}}(v,\tau)
	\end{split}
	\end{equation*}
	On the other hand, 
	$Z_{M^i_\l}(v, -1/\tau)=\tau^{\wt[v]}\sum_{j\in J}\sum_{\mu\in \Lambda_{M^j}}S_{M^i_\l M^j_{\mu}}   Z_{M^j_{\mu}}(v,\tau)$ 
	and 	 $\{Z_{M^i_{\l}}(v,\tau)\mid i\in J, \l\in \Lambda_{M^i}\}$  are linearly independent vectors in the conformal block by Theorem \ref{MTH1}.
	This implies that 
	$$S_{M^i_\l,M^j_\mu}=\frac{c_jc_i^{-1}}{|G_{M^i}|}\sum_{k\in C_{i,j}}S_{(M^i)^s,(M^j)^s\circ k}^{(g_i,k^{-1}g_j k)}\overline{\lambda(\widehat{k^{-1}g_jk})}\mu(\widehat{kg_i^{-1}k^{-1}}).$$
	
	If $i=0$, $g_0=1,$ $G_{M^0}=G$ then $\overline{\lambda(\widehat{k^{-1}g_jk})}=\lambda(g_j^{-1})$ and
	$\mu(\overline {kg_0^{-1}k^{-1}})=\dim W_{\mu}.$
	Using Lemma \ref{SMN}  yields
	$$S_{(M^0)^s,(M^j)^s\circ k}^{(1,k^{-1}g_jk)}=S_{M^0\circ k,(M^j)^s\circ k}^{(k^{-1}1k,k^{-1}g_jk)}=S_{V,(M^j)^s}^{(1,g_j)}.$$
	Thus
	$$S_{M^0_\l,M^j_\mu}=\frac{a_j}{|G_{M^j}|}S_{V,(M^j)^s}^{(1,g_j)}\lambda(g_j^{-1})\dim W_{\mu}.$$
	The proof is complete.
\qed

The following result follows from Theorem  \ref{ZVG} immediately.
\begin{coro}\label{coro1}  Let $i,j, \lambda, \mu$ be as before. If  $G$ is abelian,
	then
	$$S_{M^i_\l,M^j_\mu}=\frac{a_ja_i^{-1}}{|G_{M^i}|}\sum_{k\in C_{i,j}}S_{(M^i)^s,(M^j)^s\circ k}^{(g_i,g_j)}\overline{\lambda(\widehat{g_j})}\mu(\widehat{g_i^{-1}})$$
	when $C_{i,j}$ is not empty, and $S_{M^i_\l,M^j_\mu}=0$ otherwise.
	\end{coro}

\end{document}